\documentclass[10 pt, journal]{IEEEtran}  
\IEEEoverridecommandlockouts                             

\usepackage{xspace,amssymb,epsfig,subfigure,syntonly}
\usepackage{epsfig,amsmath,color}
\usepackage{xspace,syntonly}
\usepackage{url}
\usepackage{algorithm,algorithmic}
\usepackage{mathtools}
\usepackage{cite}

\usepackage{xspace,amssymb,epsfig,syntonly}
\usepackage{epsfig,amsmath,color}
\usepackage{xspace,syntonly,empheq}
\usepackage{url}
\usepackage{bbm}
\usepackage{mathtools}
\usepackage{cite}
\usepackage{graphicx}
\usepackage{algorithm,algorithmic}
\usepackage{verbatim}
\usepackage{amssymb}
\usepackage{amsthm}
\usepackage{enumerate}

\newcommand{\andy}[1]{\textcolor{blue}{\textsf{Andrey: #1}}\xspace}

\newcommand{\utwi}[1]{\mbox{\boldmath $#1$}}

\newcommand{\diag}{\mathsf{diag}}

\newcommand{\rev}[1]{{\color{black} #1}}

\renewcommand{\hat}{\widehat}
\renewcommand{\tilde}{\widetilde}

\newcommand{\cD}{{\cal D}}

\newcommand{\cN}{{\cal N}}
\newcommand{\cP}{{\cal P}}

\newcommand{\cC}{{\cal C}}
\newcommand{\cE}{{\cal E}}
\newcommand{\cI}{{\cal I}}

\newcommand{\cM}{{\cal M}}

\newcommand{\cV}{{\cal V}}
\newcommand{\cX}{{\cal X}}

\newcommand{\ba}{{\bf a}}
\newcommand{\bb}{{\bf b}}
\newcommand{\bd}{{\bf d}}
\newcommand{\be}{{\bf e}}

\newcommand{\bg}{{\bf g}}

\newcommand{\bp}{{\bf p}}
\newcommand{\bq}{{\bf q}}

\newcommand{\bs}{{\bf s}}
\newcommand{\bx}{{\bf x}}

\newcommand{\bv}{{\bf v}}

\newcommand{\bi}{{\bf i}}
\newcommand{\bz}{{\bf z}}
\newcommand{\by}{{\bf y}}
\newcommand{\bA}{{\bf A}}
\newcommand{\bB}{{\bf B}}

\newcommand{\bH}{{\bf H}}

\newcommand{\bM}{{\bf M}}

\newcommand{\bI}{{\bf I}}
\newcommand{\bX}{{\bf X}}

\newcommand{\bY}{{\bf Y}}

\newcommand{\bell}{{\utwi{\ell}}}

\newcommand{\bepsilon}{{\utwi{\epsilon}}}
\newcommand{\bgamma}{{\utwi{\gamma}}}

\newcommand{\blambda}{{\utwi{\lambda}}}

\newcommand{\bnu}{{\utwi{\nu}}}

\newcommand{\bPhi}{{\utwi{\Phi}}}

\newcommand{\bxi}{{\utwi{\xi}}}

\newcommand{\bmu}{{\utwi{\mu}}}
\newcommand{\bzeta}{{\utwi{\zeta}}}

\newcommand{\reals}{\mathbb{R}}
\newcommand{\comps}{\mathbb{C}}
\newcommand{\proj}{\mathrm{proj}}

\newcommand{\sfH}{\textsf{H}}
\newcommand{\sfT}{\textsf{T}}

\newcommand{\conv}{\mathsf{ch}}

\DeclarePairedDelimiterX{\norm}[1]{\lVert}{\rVert}{#1}

\usepackage{epstopdf}

\newtheorem{proposition}{Proposition}
\newtheorem{lemma}{Lemma}
\newtheorem{theorem}{Theorem}
\newtheorem{corollary}{Corollary}
\theoremstyle{definition}

\newtheorem{assumption}{Assumption}

\begin{document}

\IEEEoverridecommandlockouts

\title{Real-Time Feedback-Based Optimization of Distribution Grids:  A Unified Approach} 

\author{Andrey Bernstein and Emiliano Dall'Anese$^\star$
\thanks{$^\star$Alphabetical order, authors contributed equally to the paper. A. Bernstein is with the National Renewable Energy Laboratory, Golden, CO; e-mail: andrey.bernstein@nrel.gov. E. Dall'Anese is with the Department of Electrical, Computer, and Energy Engineering at the University of Colorado Boulder; e-mail: emiliano.dallanese@colorado.edu.}
}

\maketitle
\begin{abstract}
This paper develops an algorithmic framework for real-time optimization of distribution-level distributed energy resources (DERs). The proposed framework optimizes the operation of
\rev{both DERs that are individually controllable and groups of
DERs (i.e., aggregations) that are jointly controlled at an electrical
point of connection}. From an electrical standpoint, wye and delta single- and multi-phase connections are accounted for. The algorithm enables (groups of) DERs to pursue given performance objectives, while adjusting their (aggregate) powers to respond to services requested by grid operators and to maintain electrical quantities within engineering limits. The design of the algorithm leverages a  time-varying bi-level problem formulation capturing various performance objectives and engineering constraints, and an online implementation of primal-dual projected-gradient methods. The gradient steps are suitably modified to accommodate appropriate measurements from the distribution network and the DERs. By virtue of this approach, the resultant algorithm can cope with inaccuracies in the distribution-system modeling, it avoids pervasive metering to gather the state of non-controllable resources, and it naturally lends itself to a distributed implementation. Analytical stability and convergence claims are established in terms of tracking of the solution of the formulated time-varying  optimization problem. The proposed method is tested in a realistic distribution system with real data. 
\end{abstract}

\vspace*{-0.5cm}

\section{Introduction}
\label{sec:introduction}

This paper focuses on \emph{real-time optimization} of heterogeneous distributed energy resources (DERs) in utility-level systems and ``soft'' microgrids, with the latter referring to community-, campus-, and neighborhood-level systems connected to the rest of the grid through one point of interconnection. This paper seeks contributions in the design of real-time optimization strategies, to offer decision making capabilities that match the time scale of distribution grids with high DER integration. The objective is to allow the maximization of given DER-level and system-level operational objectives, while coping with the variability of ambient conditions and non-controllable energy assets~\cite{DhopleNoFuel}. 

Centralized and distributed optimization approaches -- such as the AC optimal power flow (OPF) -- have been developed  for distribution grids to compute optimal setpoints for DERs, so that power losses and voltage deviations are minimized and economic benefits to utility and end-users are maximized (see the tutorial~\cite{molzahn2017survey} and pertinent references therein). Centralized approaches utilize off-the-shelf solvers for nonlinear programs, or leverage convex relaxation and approximation techniques to obtain convex surrogates. Distributed solution methods capitalize on the decomposability of Lagrangian functions to decompose the solution of the optimization task across DERs, utility, and possibly ``aggregators.'' Either way, these approaches are inadequate for real-time optimization for the following main reasons: 

\noindent \emph{c1)} Computational complexity  may render \rev{impossible} the solution of optimization problems on a second (or a few seconds) timescale~\cite{molzahn2017survey,PaudyalyISGT}. In distributed settings, multiple communication rounds are required to reach convergence to a solution.  

\noindent \emph{c2)} Conventional optimization tasks operate in an \emph{open-loop} (i.e., feed-forward) setting, where a grid model and measurements of uncontrollable assets are utilized as inputs. Approximate representation of system physics, modeling errors, and uncertainty in the measurements/forecasts lead to solutions that may be in fact infeasible for the physical power system.  

\noindent \emph{c3)} Feed-forward techniques require measurements \rev{(or estimation)} of the state of non-controllable energy assets \emph{everywhere} (they are inputs of the optimization problem to be solved). Pervasive metering is impractical in existing distribution grids, \rev{and estimation based on limited data might be  inaccurate.}

\begin{figure}
  \centering
  \includegraphics[width=1.0\columnwidth]{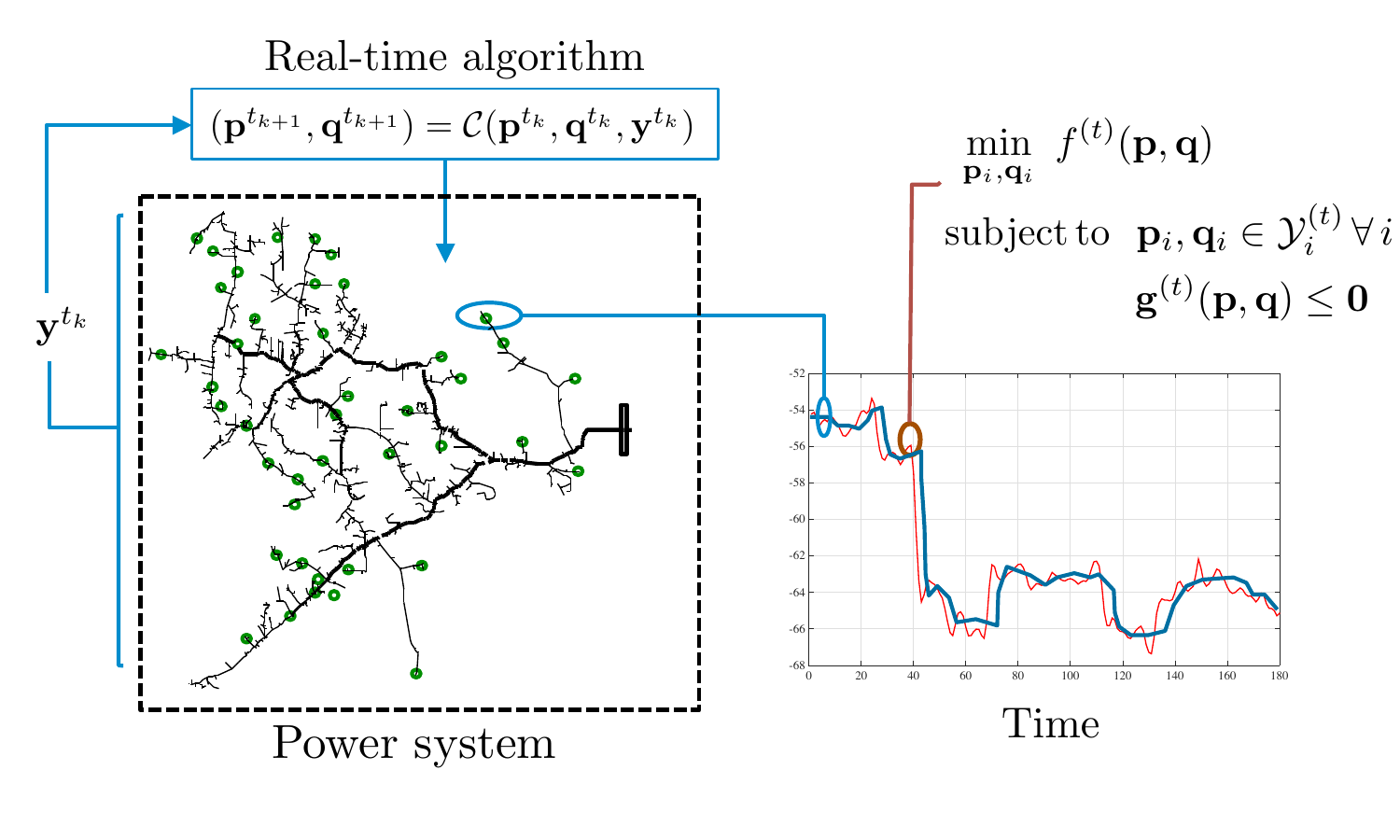}
\vspace{-.8cm}
\caption{Operating principles of the feedback-based online algorithm. Power setpoints $\bp, \bq$ of the devices are updated in real-time through the map $\cC$. The design of the update $\cC(\bp, \bq, \by)$ capitalizes on online implementations of gradient-based methods, suitably modified to accommodated feedback $\by$ (i.e., measurements) from the power system and the devices. Analytical convergence claims demonstrate the tracking of the solution of a time-varying  optimization problem.} \label{fig:F_feedback}
\vspace{-.5cm}
\end{figure}

This paper starts from the formulation of a \emph{time-varying bi-level convex optimization problem} that models optimal operational trajectories of DERs and groups of DERs, and embeds dynamic operational and engineering constraints. The latter include voltage and ampacity limits, feasible operating regions of DERs, and target power flows at the point of interconnection with the rest of the grid (to provide services~\cite{Caiso12} or partake into market operations). To address the challenge c1), we develop an online algorithm based on a projected primal-dual gradient method to track the optimal solution of the formulated optimization problem. To resolve c2) and c3), the gradient steps are suitably modified to accommodate voltage, power, and current measurements from the distribution network and the DERs -- hence the term \emph{feedback-based online optimization}. The operating principles of the real-time framework are illustrated in Fig.~\ref{fig:F_feedback}. 
The synthesis of the algorithm leverages the fixed-point linearization method for the multi-phase AC power-flow equations presented in~\cite{multiphaseArxiv}, where delta and wye connections are unified under the same mathematical formalism. The resultant algorithm avoids pervasive metering to gather the state of non-controllable resources, it can cope with inaccuracies in the representation of the AC power flows, and it affords a distributed implementation. As illustrated  in Fig.~\ref{fig:F_feedback}, analytical convergence and stability claims are established in terms of tracking of the solution of the formulated time-varying  optimization problem. 

The bi-level nature of the problem allows one to readily optimize the \emph{net} power generated/consumed by groups of DERs located behind the same meter, while accounting for individual DER constraints. Towards this end, the paper contributes results with respect to the computation of inner approximations of the Minkowski sum of prototypical sets of DERs, to represent the overall feasibility region of groups of DERs; and, a systematic way to compute the gradient of the cost function associated with groups of DERs, along with a  mechanism to disaggregate the power command across DERs. With respect to the types of DERs, the paper considers DERs with both continuous and discrete implementable power commands. For the latter, the operational sets of DERs are convexified for the purpose of setpoint computation, whereas implementable setpoints are computed via error-diffusion techniques \cite{Anastassiou,bernstCDC}. 

\rev{General bi-level optimization problems are  closely related to the Stackelberg game formulation \cite{von1952theory}, and are NP-hard even when both the inner and outer optimization problems are convex \cite{Bard1991}. There is an extensive literature on such problems in different application domains \cite{Colson2007,Tseveendorj2013}. In our application, the solution to the inner level optimization problem can be expressed as the solution to the corresponding KKT first order optimality conditions, and hence the problem can be efficiently solved.}

The idea of leveraging \emph{time-varying problem formulations} to model optimal operational trajectories for DERs and developing feedback-based online solvers to track the optimal trajectory traces back to our preliminary works~\cite{commelec1,commelec2} and~\cite{opfPursuit}, where time-varying linearized AC OPF formulations were considered for distribution grids. Feedback was in the form of measurements of voltages and powers; while the effectiveness of these methods were shown numerically in~\cite{commelec1,commelec2}, analytical tracking results were first provided in~\cite{opfPursuit}. A centralized online algorithm based on a quasi-Newton method was proposed in~\cite{Tang17} for a time-varying relaxed non-convex AC OPF (smooth penalty functions were utilized to relax the constraints); considerations regarding estimations of the Hessian were offered.  An online incentive-based algorithm was developed in~\cite{Zhou17} to track a time-varying equilibrium point of a Stackelberg game. Voltage measurements were included into the steps of the alternating direction method of multipliers in~\cite{ZhangADMM17}, and tracking results were established. 

For \emph{static} optimization problems, a feedback-based algorithm for a real-time solution of economic dispatch problems was proposed in~\cite{Jokic_JEPES};  feedback was in the form of measurements of the output powers. Measurements of voltages were considered in the distributed strategy developed in~\cite{Christakou14}, to enforce voltage regulation in distribution networks. Similarly,~\cite{Bolognani_feedback_15} proposed a distributed  \rev{reactive power control strategy, and convergence} to a solution of a well-posed static optimization problem was analytically established. Local control methods for voltage regulation were proposed in~\cite{Zhu16_local_control}; performance in a dynamic setting was experimentally evaluated. Power measurements were utilized in~\cite{DhopleDKKT15} to dynamically solve a relaxed AC OPF. State measurements were leveraged in~\cite{LowOnlineOPF} to solve an AC OPF for radial systems. Manifold-based approaches were proposed and analyzed in~\cite{Hauswirth16,Hauswirth17} to solve the AC OPF (smooth penalty functions were utilized to relax the constraints); 
however, the update of the tangent plane in these papers may still require pervasive metering of the non-controllable assets.

The framework in this paper significantly expands our prior works in~\cite{opfPursuit,optimalregulationvpp2017} by providing the following  contributions:  

\rev{
\noindent \emph{i)} We consider a bi-level optimization formulation, which is NP-hard problem in general. In particular, the rigorous analysis of \emph{online bi-level formulation with feedback} is absent in the literature.

\noindent \emph{ii)} We account for aggregations of DERs (e.g., buildings and facilities with multiple DERs behind the meter). We provide new results for the inner approximation of the Minkowski sum of prototypical operational sets for DERs, to represent the overall feasibility set of groups of DERs; we also offer new insights on the computation of the gradient of the aggregate cost function associated with groups of DERs. These results are of independent interest for real-time control applications.      

\noindent \emph{iii)} The proposed algorithm is applicable to multi-phase systems with both wye and delta connections.

\noindent \emph{iv)} The proposed framework accommodates DERs with a non-convex (and, in particular, discrete) set of implementable control commands.

\noindent \emph{v)} The proposed framework is tested through numerical simulations based on a real system and using real data from a distribution network located in California in the territory of Southern California Edison. 
}


\section{Preliminaries and System Model}
\label{sec:model}

\emph{Notation}: Upper-case (lower-case) boldface letters will be used for matrices (column vectors); $(\cdot)^\sfT$ for transposition; $(\cdot)^*$ complex-conjugate; and, $(\cdot)^\sfH$ complex-conjugate transposition. $\Re\{\cdot\}$ and $\Im\{\cdot\}$ denote the real and imaginary parts of a complex number, respectively, and $\mathrm{j} := \sqrt{-1}$. $|\cdot|$ denotes the absolute value of a number or the cardinality of a (discrete) set.  For a given $N \times 1$ vector $\bx \in \mathbb{R}^N$, $|\bx|$ takes the absolute value entry-wise; $\|\bx\|_2 := \sqrt{\bx^\sfH \bx}$; and, $\diag(\bx)$ returns a $N \times N$ matrix with the elements of $\bx$ in its diagonal.  Given a matrix $\bX \in \mathbb{R}^{N\times M}$, $x_{m,n}$  denotes its $(m,n)$-th entry and $\|\bX\|_2$ denotes the $\ell_2$-induced matrix norm. For a function $f: \mathbb{R}^N \rightarrow \mathbb{R}$, $\nabla_{\bx} f(\bx)$ returns the gradient vector of $f(\bx)$ with respect to $\bx \in \mathbb{R}^N$. The symbols $\mathbf{1}_N$ and $\mathbf{0}_N$ denote the $N \times 1$ vector with all ones and with all zeros, respectively. Given two sets $\cX_1 \subset \mathbb{R}^N$ and $\cX_2 \subset \mathbb{R}^N$, $\cX_1  \oplus  \cX_2$ denotes the Minkowski sum of $\cX_1$ and $\cX_2$. Finally, given a set $\cX \subset \mathbb{R}^N$, $\conv \cX$ denotes its convex hull; and $\mathrm{proj}_{\cX}\{\bx\}$ denotes a closest point to $\bx$ in $\cX$, namely $\mathrm{proj}_{\cX}\{\bx\} \in \arg \min_{\by \in \cX} \|\bx - \by\|_2 $ (the ties can be broken arbitrarily).

\vspace{-0.5cm}

\subsection{DER Model}
\label{sec:modelDER}

We consider two classes of DERs: \emph{i)} devices that are individually controllable; and, \emph{ii)} groups of DERs that can be controlled as a whole. The second class models, e.g., residential homes and buildings with multiple DERs behind the meter, renewable-based systems with multiple (micro)inverters, and parking garages for EVs.  Each DER can be either wye-connected or delta-connected to the network~\cite{Kerstingbook}, and it can be either single-phase or three-phase. In the following, pertinent notation and modeling details are outlined.

For future developments, let $\cP := \{a, b, c\} \cup \{ab, bc, ca\}$ be the set of possible connections, with $\{a, b, c\}$ pertaining to wye connections (line to ground) and $\{ab, bc, ca\}$ referring to delta connections (line to line). 

\textbf{\emph{Individually-controllable DERs}}. Let $\cD := \{1,\ldots, D\}$ be the set of individually-controllable DERs, and let \noindent $\bx_j := [P_j, Q_j]^\sfT \in \mathbb{R}^2$ collect the real and reactive power setpoint of  DER $j \in \cD$. The DER can be either wye-connected or delta-connected to the network. Three-phase DERs are assumed to operate in a balanced setting; thus, the setpoint $\bx_j$ is the same across phases. The set $\cP_j \subset \cP$ collects the phases where DER $j$ is connected. 

We denote as $\cX_j \subset \reals^2$ the set of possible power setpoints $\bx_j$ for the DER $j$; the set $\cX_j$ captures hardware and operational constraints and it is assumed to be convex and compact.  It is assumed that the DERs are endowed with  controllers that are designed so that, upon receiving the setpoint $\bx_j \in \cX_j$, the output powers are driven to the commanded setpoints; relevant dynamical models for the output powers of an inverter-interfaced DER are discussed in e.g.,~\cite{Iravanibook10,Irminger12} and can be found in  datasheets of commercially available DERs. 

For an inverter-interfaced DER, we consider the following prototypical representation of the set $\cX_j$: 
\begin{align} 
\label{eq:set}
\hspace{-.2cm} \cX_j(\underline{p},\overline{p},r)  := \left\{ [P_j, Q_j]^\sfT \hspace{-.1cm} : \underline{p} \leq P_j  \leq  \overline{p}, P_j^2 + Q_j^2 \leq  r^2 \right\}
\end{align}
where $\underline{p}$, $\overline{p}$, and $r > 0$ are given DER-dependent parameters. For example,  for a PV system $r$ represents the inverter capacity, $\underline{p} = 0$, and $\overline{p}$ is the available real power. For an energy storage systems, $r$ represents the inverter rating, and $\underline{p}, \overline{p}$ are updated during the operation of the battery based on its internal state (such as the  state of charge or DC voltage). Notice that the set $\cX_j$ is typically \emph{time varying}, as the parameters $\underline{p}$, $\overline{p}$, and $r$ vary over time based on ambient conditions and/or internal DER state. 

On the other hand, we consider the following  operating region for DERs with controllable active powers   (e.g.,  variable speed drives and EVs; see other examples in~\cite{GatsisTSG12}): 
\begin{align} 
\label{eq:setP}
\cX_j(\underline{p},\overline{p})  := \left\{ [P_j, Q_j]^\sfT \hspace{-.1cm} :  \underline{p} \leq P_j  \leq  \overline{p},  Q_j = 0\right\}.
\end{align}

\textbf{\emph{DERs with nonconvex (discrete) controls}}. Consider a DER with a nonconvex operating region, $\tilde{\cX}_j \subset \reals^2$. This is the case, for example, for \rev{(residential)} HVAC systems where $\tilde{\cX}_j = \{[P_j, Q_j]^\sfT:  P_j \in \{0, \overline{p}\}, Q_j = 0\}$, or EVs with discrete charging levels. For these devices, the set $\cX_j$ is \emph{the convex hull} of $\tilde{\cX}_j$; i.e., $\cX_j := \textrm{ch}\tilde{\cX}_j$. For example,  for an HVAC systems we have that $\cX_j  =  \{[P_j, Q_j]^\sfT: 0\leq  P_j \leq \overline{p}, Q_j = 0\}$ [cf.~\eqref{eq:setP}].  The algorithm proposed in Section~\ref{sec:algorithm} will utilize a randomization procedure to recover implementable setpoints based on $\cX_j$~\cite{floyd75,errDiff}.  For a DER with nonconvex set of implementable setpoints, $\tilde{\bx}_j \in \tilde{\cX}_j$ denotes an implementable setpoint, whereas 
$\bx_j \in \cX_j$ is a (relaxed) setpoint computed based on the convex hull of $\tilde{\cX}_j$. Notice that for devices that lock on a state due to engineering or operational constraints, the set $\tilde{\cX}_j$ is a singleton over a given period of time; for example, if an HVAC system is required not to switch ON for a few minutes, then $\tilde{\cX}_j(t) = \{0\}$ for a given interval $t \in [t_0, t_0 + T]$.

\textbf{\emph{Aggregations of DERs}}. Suppose that the distribution grid features a set $\bar{\cD} := \{1,\ldots, \bar{D}\}$ of residential  homes, building, or other facilities with multiple DERs that are controlled jointly. Let $\bar{\cD}_j := \{1, \ldots, \bar{D}_j\}$ denote the set of devices within the $j$th aggregation, and define as $\bar{\bx}_j := \sum_{i \in \bar{\cD}_j} \bx_i$ the setpoint for the \emph{net powers generated} by the DERs within a group. The set $\bar{\cP}_j \subset \cP$ collects the connections of the aggregation $j$. 

Let $\bar{\cX}_j \subseteq \oplus_{i \in \bar{\cD}_j} \cX_i$ be either the exact Minkowski sum of the operating regions of the DERs within the $j$th aggregation or an inner approximation thereof.  
In the following, we provide pertinent results for the Minkowski sum of sets~\eqref{eq:set} and~\eqref{eq:setP}. \rev{There is extensive literature on devising \emph{generic} numerical methods to compute aggregated flexibility regions using Minkowski sums \cite{Muller2015,Kundu2018,Nazir2018}. The goal in this paper is rather to develop a \emph{simple analytical approximation} that can be directly used in \emph{real-time algorithms}.} 

First, notice that the  Minkowski sum of two sets $ \cX_j(\underline{p}_j,\overline{p}_j)$ and $\cX_n(\underline{p}_n,\overline{p}_n)$ for two DERs with controllable active powers is given by:
\begin{align}
\label{eq:sumP}
& \cX_j(\underline{p}_j,\overline{p}_j) \oplus \cX_n(\underline{p}_n,\overline{p}_n) \nonumber \\
& \hspace{.5cm} = \left\{[P, Q]^\sfT \hspace{-.1cm} :  \underline{p}_j + \underline{p}_n  \leq P  \leq  \overline{p}_j + \overline{p}_n,  Q = 0\right\}.
\end{align}
The following theorems deal with the Minkowski sums  $\cX_j(\underline{p}_j,\overline{p}_j, r_j) \oplus \cX_n(\underline{p}_n,\overline{p}_n)$ and $\cX_j(\underline{p}_j,\overline{p}_j, r_j) \oplus \cX_n(\underline{p}_n,\overline{p}_n, r_n)$. 

\rev{
\begin{proposition}
\label{thm:sum1}
The Minkowski sum between $\cX(\underline{p}_1,\overline{p}_1, r)$ and $\cX(\underline{p}_2,\overline{p}_2)$ in~\eqref{eq:set} and~\eqref{eq:setP}, respectively, with $\underline{p}_1 \in [-r, 0], \overline{p}_1 \in [0, r]$, is given by
\begin{align}
& \hspace{-.5cm} \cX(\underline{p}_1,\overline{p}_1, r) \oplus \cX(\underline{p}_2,\overline{p}_2) = \Big \{[P, Q]^\sfT: \underline{p}_1+\underline{p}_2 \leq P \leq \overline{p}_1 + \overline{p}_2,  \nonumber\\
& \hspace{1.3cm} \qquad -g(P) \leq Q \leq g(P) \Big\},
\label{eq:sum1}
\end{align}
where $g(P)$ is a concave function given by:
\[
g(P) :=
\begin{cases}
r, &  P \in [\underline{p}_2, \overline{p}_2], \\
\sqrt{r^2 - (P - \underline{p}_2)^2 }, & P \in [\underline{p}_1+\underline{p}_2, \underline{p}_2), \\
\sqrt{r^2 - (P - \overline{p}_2)^2 }, & P \in (\overline{p}_2, \overline{p}_1 + \overline{p}_2];\\
\end{cases}
\] 
we use the convention that an interval $[a, b)$ (or $(a, b]$) with $a > b$ is an empty set.
\hfill $\Box$
\end{proposition}
}

\begin{proposition}
\label{thm:sum2}
Inner and outer approximations of the Minkowski sum of two sets $\cX(\underline{p}_1,\overline{p}_1, r_1)$ and $\cX(\underline{p}_2,\overline{p}_2, r_2)$ are given by
\begin{subequations}
\begin{align}
& \hspace{-.2cm} \cX(\underline{p}_1 + \underline{p}_2, \overline{p}_1 + \overline{p}_2, \rho)  \subseteq 
\cX(\underline{p}_1,\overline{p}_1, r_1) \oplus \cX(\underline{p}_2,\overline{p}_2, r_2) \label{eq:sum_inner} \\
& \hspace{3.05cm} \subseteq 
\cX(\underline{p}_1 + \underline{p}_2, \overline{p}_1 + \overline{p}_2, r_1 + r_2) \label{eq:sum_outer} 
\end{align}
\end{subequations}
for any $\rho > 0$ satisfying the following condition 
\begin{equation}
\label{eq:rho}
\rho^2 \leq r_1^2 + r_2^2 + \alpha - \beta_1 - \beta_2 + 2\sqrt{(r_1^2 - \beta_1)(r_2^2 - \beta_2)}, 
\end{equation}
where $\alpha := [ \max \{\underline{p}_1 + \underline{p}_2, \min \{0,  \overline{p}_1 + \overline{p}_2,\} \} ]^2$, and $\beta_i := \max \{\underline{p}_i^2, \overline{p}_i^2 \}, \quad i = 1, 2$. \hfill $\Box$
\end{proposition}%
\rev{We note that the best choice for $\rho$ in Proposition \ref{thm:sum2} is given by the upper bound:
\[
\rho = \sqrt{r_1^2 + r_2^2 + \alpha - \beta_1 - \beta_2 + 2\sqrt{(r_1^2 - \beta_1)(r_2^2 - \beta_2)}}.
\]}%
The proofs of \rev{Propositions}~\ref{thm:sum1} and \ref{thm:sum2} are provided in Appendix~\ref{sec:proofthm1} and \ref{sec:proofthm2}, respectively.
Notice that the inner approximation $\cX(\underline{p}_1 + \underline{p}_2, \overline{p}_1 + \overline{p}_2, \rho) $ is convex and compact, and it is in the form of~\eqref{eq:set}. 

Expression~\eqref{eq:sumP} along with the results of Propositions~\ref{thm:sum1} and Theorem~\ref{thm:sum2} can be utilized to compute an inner approximation of the feasible region of the net powers $\bar{\bx}_j$ for each aggregation of DERs $j \in \bar{\cD}$. For example, the feasible region for the net powers generated by a residential house with a PV system, a battery, and an EV can be computed by first leveraging~\eqref{eq:sum_inner} to sum the sets pertaining to the PV system and the battery and subsequently~\eqref{eq:sum1}, to add up the feasible region of the EV.   

\vspace*{-0.5cm}

\subsection{Network Model}
\label{sec:modelNetwork}

We consider a generic multi-phase distribution network with multiphase nodes collected in the set $\cN \cup \{0\}$, $\cN := \{1,\ldots, N\}$, and distribution lines represented by the set of edges $\cE := \{(m,n)\} \subset (\cN  \cup \{0\}) \times (\cN  \cup \{0\})$. Node $0$ denotes the three-phase slack bus, i.e., the point of connection of the distribution grid with the rest of the electrical system. At each multiphase node, controllable and non-controllable devices can be either wye- or delta-connected~\cite{Kerstingbook}. 

We briefly showcase the set of AC power-flow equations for this generic setting (a comprehensive description can be found in, e.g.,~\cite{Kerstingbook} and~\cite{multiphaseArxiv}). To this end, let $\bv$ be a vector collecting the line-to-ground voltages in all phases of the nodes in $\cN$; similarly, vector $\bi$ collects all the phase net current injections, $\bi^{\Delta}$ the phase-to-phase currents in all the delta connections, and vectors $\bs^Y$ and $\bs^\Delta$ collect the net complex powers injected at nodes $\cN$ from devices with wye and delta connections, respectively. With these definitions in place, the AC power-flow equations can be compactly written as:   
\begin{subequations} \label{eq:lf}
\begin{align} 
\diag(\bH^\sfT (\bi^{\Delta})^*)\bv + \bs^Y = \diag(\bv) \bi^*, \label{eq:lf_balance}\\
\bs^{\Delta} = \diag\left(\bH \bv \right)(\bi^{\Delta})^*, \, \, 
\bi = \bY_{L0} \bv_0 + \bY_{LL} \bv , \label{eq:lf_i} 
\end{align}
\end{subequations}
where $\bY_{00} \in \comps^{3 \times 3}, \bY_{L0} \in \comps^{N_\phi \times 3}, \bY_{0L} \in \comps^{3 \times N_\phi}$, and $\bY_{LL} \in  \comps^{N_\phi \times N_\phi}$ are the submatrices of the admittance matrix
\begin{equation}
\bY := 
\begin{bmatrix}
\bY_{00} & \bY_{0L} \\
\bY_{L0} & \bY_{LL}
\end{bmatrix} \in \comps^{N_\phi+3 \times N_\phi + 3},
\end{equation}
which can be formed from the topology of the network and the $\pi$-model of the distribution lines~\cite{Kerstingbook}; $N_\phi$ is the total number of single-phase connections, and $\bH$ is a $N_\phi \times N_\phi$ block-diagonal matrix mapping the direction of the currents $\bi^{\Delta}$ in the delta connections; see~\cite{multiphaseArxiv,ZhaoIREP17} for a detailed description. 

The nonlinearities in~\eqref{eq:lf} hinder the possibility of seeking analytical closed-form solutions to pertinent electrical quantities such as voltages, power flows, and line currents as a function of the DERs' power injections. To facilitate the design and  analysis of real-time optimization methods, we leverage the approximate linear models proposed in~\cite{multiphaseArxiv,linModels}. To this end, denote as $\bv_{\cM_v}$ the vector collecting the phase-to-ground voltages at given measurement points $\cM_v$; $\bi_{L,\cM_i}$ the vector  collecting the line currents for a subset of monitored distribution lines $\cM_i$ (or given by pseudo-measurements); and, $\bp_0 \in \reals^3$ the vector of real powers entering node $0$ on the phases $\{a, b, c\}$. Then, $|\bv_{\cM_v}|$ (where the absolute value is taken entry-wise), $|\bi_{L,\cM_i}|$, and $\bp_0$ can be approximately expressed as:  
\begin{subequations} 
\label{eq:approx_v}
\begin{align}
|\tilde{\bv}_{\cM_v}(\bx, \bar{\bx})| & = \sum_{j \in \cD}  \bA_{j} \bx_j + \sum_{j \in \bar{\cD}}  \bar{\bA}_{j} \bar{\bx}_j + \ba 
\label{eq:approx_v1} \\
\mathbf{a} & := \sum_{j \in \cN} \sum_{\phi \in \cP} \bA_{j,\phi} \bell_{j,\phi} + \mathbf{a}_0 \label{eq:approx_v2}
\end{align}
\end{subequations}
\vspace*{-0.2cm}
\begin{subequations} 
\label{eq:approx_i}
\begin{align}
|\tilde{\bi}_{L,\cM_i}(\bx, \bar{\bx})| & = \sum_{j \in \cD}  \bB_{j} \bx_j + \sum_{j \in \bar{\cD}}  \bar{\bB}_{j} \bar{\bx}_j + \bb 
\label{eq:approx_v1} \\
\mathbf{b} & := \sum_{j \in \cN} \sum_{\phi \in \cP} \bB_{j,\phi} \bell_{j,\phi} + \mathbf{b}_0 \label{eq:approx_v2}
\end{align}
\end{subequations}
\begin{subequations} 
\label{eq:approx_p}
\begin{align}
\tilde{\bp}_0(\bx, \bar{\bx}) & =  \sum_{j \in \cD} \bM_{j} \bx_j + \sum_{j \in \bar{\cD}}  \bar{\bM}_{j} \bar{\bx}_j + \mathbf{m} 
\label{eq:approx_p01} \\
\mathbf{m} & := \sum_{j \in \cN} \sum_{\phi \in \cP} \bM_{j,\phi} \bell_{j,\phi} + \mathbf{m}_0 \label{eq:approx_p02}
\end{align}
\end{subequations}
where   $\bell_{j,\phi} \in \reals^2$ collects the net non-controllable active and reactive powers at connection $\phi \in \cP$ of node $n \in \cN$, $\bx$ and $\bar{\bx}$ stack all the setpoints $\{\bx_j\}$ and $\bar{\bx}_j$, respectively, and the matrices $\bA_{j,\phi}$, $\bar{\bA}_{j,\phi}$, $\bB_{j,\phi}$, $\bar{\bB}_{j,\phi}$, $\bM_{j,\phi}$, $\bar{\bM}_{j,\phi}$ along with the vectors $\mathbf{a}_0$, $\mathbf{b}_0$, and $\mathbf{m}_0$ are  model parameters that can be computed through e.g., the fixed-point linearization  method  proposed in~\cite{multiphaseArxiv,linModels}; for brevity, we defined the matrices $\bA_{j} := \sum_{\phi \in \cP_j} \bA_{j,\phi}$, $\bar{\bA}_{j} := \sum_{\phi \in \bar{\cP}_j} \bar{\bA}_{j,\phi}$, $\bB_{j} := \sum_{\phi \in \cP_j} \bB_{j,\phi}$, $\bar{\bB}_{j} := \sum_{\phi \in \bar{\cP}_j} \bar{\bB}_{j,\phi}$, $\bM_{j} := \sum_{\phi \in \cP_j} \bM_{j,\phi}$, and $\bar{\bM}_{j} := \sum_{\phi \in \bar{\cP}_j} \bar{\bM}_{j,\phi}$. As explained in~\cite{multiphaseArxiv,linModels}, these model parameters capture the effects of different types of connection (e.g., wye or delta) and can be computed based on the admittance matrix of the system. If a  fixed-point linearization  method is utilized, knowledge of the non-controllable powers $\bell_{j,\phi}$ is not required for the computation of the  model parameters. If only wye connections are present, an alternative way to obtain~\eqref{eq:approx_v}--\eqref{eq:approx_p} is presented in, e.g., \cite{christ2013sens}. 

It is worth emphasizing that the approximate models~\eqref{eq:approx_v}--\eqref{eq:approx_p} are utilized to facilitate the design and the performance analysis of the real-time algorithm. 
In Section~\ref{sec:algorithm}, we  show how to leverage  measurements from the distribution grid and DERs to cope with the inaccuracies \rev{introduced by a linear approximation of the AC power flows; whereas in Section \ref{sec:convergence}, we establish appropriate stability claims}. 

Hereafter, we will drop the subscripts $\cM_v$ and $\cM_i$ from~\eqref{eq:approx_v} and~\eqref{eq:approx_i} for notational simplicity, with the understanding that functions $\bv(\bx, \bar{\bx})$ and $\bi_L(\bx, \bar{\bx})$ refer to voltages and currents at given points of interest. 

\vspace*{-0.2cm}

\section{Feedback-based Dynamic Optimization}
\label{sec:algorithm}

We design a new \emph{real-time optimal power flow} method where power setpoints of the DERs are updated on a second timescale~\cite{AndreayOnlineOpt, opfPursuit,Tang17} to  maximize operational objectives while coping with the variability of ambient conditions and non-controllable assets. Consider then discretizing the temporal domain as $t_k = k h$, where $k \in \mathbb{N}$ and $h > 0$ will be taken to be the time required to compute one closed-loop iteration of the proposed algorithm. As discussed shortly, the value of $h$ is based on underlying communication  delays, as well as operational considerations of utility and aggregators.  

We next leverage the \emph{time-varying optimization} formalism~\cite{SimonettoGlobalsip2014,opfPursuit} to model optimal operational trajectories for the DERs, based on 1) possibly time-varying optimization objectives and operational constraints, as well as 2) variability of non-controllable assets and ambient conditions. Hereafter, the superscript $^{(k)}$ will be  utilized to indicate variables, functions, and inputs at time $t_k$, for all $k \in \mathbb{N}$. 

\vspace*{-0.5cm}

\subsection{Formalizing Optimal Operational Trajectories}

Let $v^{min}$ and $v^{max}$ be given limits for the magnitude of phase-to-ground voltages  (e.g., ANSI C.84.1 limits), and let $ \bi^{max}$ be a vector collecting the ampacity limits for the monitored distribution lines. Finally, $s^{(k)} \in \{0,1\}$ indicates whether the distribution grid is requested to  follow a setpoint $\bp_{0,{\mathrm{set}}}^{(k)}$ for the real powers at the three phases of the point of connection with the rest of the electrical network~\cite{Caiso12,Meliopoulos16}. When $s^{(k)} = 1$, the sequence of setpoints $\{\bp_{0,{\mathrm{set}}}^{(k)}\}_{k}$ shall be tracked within a given accuracy $E^{(k)}$. With these definition, the following time-varying optimization problem is formulated to model optimal operational trajectories $\{\bx_j^{\textrm{opt}}, k \in \mathbb{N}\}$ for the DERs: 
\vspace*{-0.2cm}
\begin{subequations} 
\label{Pmg2}
\begin{align} 
 \mathrm{(P1}^{(k)} \mathrm{)}  \hspace{1.8cm} & \hspace{-1.7cm} \min_{\bx,  \bar{\bx} } \,\, \sum_{j \in \cD} f_j^{(k)}(\bx_j) + \sum_{j \in \bar{\cD}} \bar{f}_j^{(k)}(\bar{\bx}_j) \label{mg-cost2} \\ 
& \hspace{-2.5cm} \mathrm{subject\,to:}~   \bx_j \in  \cX_j^{(k)} \hspace{1.5cm} \forall \, j \in \cD\   \label{mg-PVp2}   \\ 
&  \hspace{-0.7cm} \bar{\bx}_j \in  \bar{\cX}_j^{(k)} \hspace{1.5cm} \forall \, j \in \bar{\cD}\\\
& \hspace{-0.65cm}  s^{(k)} \bI_3 (\tilde{\bp}_0^{(k)}(\bx,  \bar{\bx}) - \bp_{0,{\mathrm{set}}}^{(k)})   \leq  E^{(k)} \mathbf{1}_3  \label{mg-power-sub2}  \\
& \hspace{-0.65cm}  s^{(k)} \bI_3 (\bp_{0,{\mathrm{set}}}^{(k)} - \tilde{\bp}_0^{(k)}(\bx,  \bar{\bx}))  \leq  E^{(k)} \mathbf{1}_3  \label{mg-power-sub3} \\
& \hspace{-0.65cm} |\tilde{\bv}^{(k)}(\bx, \bar{\bx})|  \leq v^{max} \mathbf{1}  \label{mg-volt1} \\
& \hspace{-0.65cm} v^{min} \mathbf{1} \leq |\tilde{\bv}^{(k)}(\bx, \bar{\bx})| \label{mg-volt2} \\
& \hspace{-0.65cm} |\tilde{\bi}_L^{(k)}(\bx, \bar{\bx})| \leq \bi^{max} \label{mg-current}
\end{align}
\end{subequations}
where we recall that $\cX_j^{(k)}$ is a convex set modeling hardware constraints of the DER $j$ at a given time $t_k$; $f_j^{(k)}:  \reals^2 \rightarrow \reals$ is a time-varying convex function associated with the  DER $j \in \cD$;  and, the function $\bar{f}_j^{(k)}: \reals^2 \rightarrow \reals$ associated with the $j$th aggregation of DERs is defined as follows:
\begin{subequations} 
\label{eq:aggregation}
\begin{align} 
\bar{f}_j^{(k)}(\bar{\bx}_j)  := \quad  &\min_{\{\bx_i\}_{i \in \bar{\cD}_j}}  \hspace{.2cm} \sum_{i \in \bar{\cD}_j} f_i^{(k)}(\bx_i)  \\
&\hspace{-.3cm} \textrm{subject to:} ~ \bx_{i} \in \cX_{i}^{(k)}, \, \forall \,\, i \in \bar{\cD}_j \label{eqn:constr_Y} \\
&\hspace{1.2cm}  \sum_{i \in \bar{\cD}_j} \bx_i = \bar{\bx}_j  \, .\label{eqn:constr_sum}
 \end{align}
\end{subequations}
Problem~\eqref{eq:aggregation} is utilized to disaggregate the setpoint $\bar{\bx}_j$ across the DERs $i \in \bar{\cD}_j$.  

Before proceeding, it is worth emphasizing the following for the bi-level formulation~\eqref{Pmg2}--\eqref{eq:aggregation}: 

\noindent i) when set $\bar{\cX}_j^{(k)}$ is given by the (exact) Minkowski sum of $ \cX_{i}^{(k)}, i \in \bar{\cD}_j$, ~\eqref{Pmg2}--\eqref{eq:aggregation} is equivalent to a ``flat'' optimization strategy where~\eqref{Pmg2} does not consider points of aggregation (the flat formulation  includes individual optimization variables and constraints for each one of the DERs; see e.g.,~\cite{LowOnlineOPF,commelec1,opfPursuit,GatsisTSG12});

\noindent ii) if the set $\bar{\cX}_j^{(k)}$ is an inner approximation of the Minkowski sum, then~\eqref{Pmg2}--\eqref{eq:aggregation} represents a restriction of the ``flat'' optimization problem. 

Problem $(\textrm{P1}^{(k)})$ is a time-varying convex optimization problem; however, solving $(\textrm{P1}^{(k)})$ in a batch fashion at each time $t_k$ might be impractical because of the following three main challenges:

\noindent $\bullet$ \textbf{c1: Complexity.} For real-time implementations (e.g., when $h$ is on the order of a second or a few seconds), it might be \rev{impossible} to solve $(\textrm{P1}^{(k)})$ to convergence; this is especially the case of distributed settings, where multiple communication rounds are \rev{required} to reach convergence.  

\noindent $\bullet$ \textbf{c2. Model inaccuracy.} 
Feed-forward (i.e., open-loop) solution of $(\textrm{P1}^{(k)})$ suffers from inaccuracies due approximate linear models
~\eqref{eq:approx_v}--\eqref{eq:approx_p}, as well as estimation errors for the admittance matrix and loads.


\noindent $\bullet$ \textbf{c3. Pervasive metering.}  Solving $(\textrm{P1}^{(k)})$ (either in a batch form or online) requires collecting measurements of the (aggregate) noncontrollable loads $\bell_{j,\phi}$ at all locations in real time, in order to compute \eqref{eq:approx_v}-\eqref{eq:approx_p}~\cite{molzahn2017survey}. 

In the following, we present a \emph{feedback-based online algorithm} that tracks the optimal solution of $\textrm{(P1)}^{(k)}$ over time, while coping with model inaccuracies and avoiding ubiquitous metering.  

\vspace*{-1cm}

\subsection{Real-time Algorithm}

The following assumption is imposed throughout the paper.

\begin{assumption} \label{ass:agg_cost}
For each DER $ i \in \cD$, and for each DER $i \in \bar{\cD}_j$ in the aggregation $j \in \bar{\cD}$:

\noindent \textbf{A\ref{ass:agg_cost}.i} the set $\cX_i^{(k)}$ is convex and compact for all $t_k$;

\noindent \textbf{A\ref{ass:agg_cost}.ii} the function $f_i^{(k)}(\bx_i)$ is convex and continuously differentiable, and its gradient is Lipschitz continuous for all $t_k$. 
\end{assumption}

\rev{
\begin{assumption} \label{ass:dual_variable}
 For problem~\eqref{eq:aggregation}, let  $d_j^{(k)}(\bxi)$ be the dual function associated with problem~\eqref{eq:aggregation} at time $t_k$~\cite{Yu18}, where $\bxi \in \reals^2$ is the Lagrange multiplier associated with constraint \eqref{eqn:constr_sum}. For any $\bar{\bx}_j^{(k)}$ in the interior of the Minkowski sum of $\cX_{i}^{(k)}, \,  i \in \bar{\cD}_j$, this dual function is locally strongly concave around an optimal dual variable $\bxi^*$.
\end{assumption}

See Appendix \ref{sec:ass2_illustration} for an example of problems that satisfy the assumption above; see also~\cite{Yu18} for analytical conditions and more elaborate examples.
}

We next outline results pertaining to the DER aggregations $\bar{\cD}$. 


\begin{lemma}
\label{lemma:strong_dual}
Suppose that problem~\eqref{eq:aggregation} is feasible and Assumptions~\ref{ass:agg_cost}-\ref{ass:dual_variable} hold. Then, the unique optimal dual variable associated with~\eqref{eqn:constr_sum} is bounded. \hfill $\Box$
\end{lemma}


\begin{theorem}
\label{prop:derivative}
Under Assumptions~\ref{ass:agg_cost}-\ref{ass:dual_variable}, it holds that:
\begin{enumerate}[(i)]
\item The function $\bar{f}_j^{(k)}(\bar{\bx})$ given in~\eqref{eq:aggregation} is convex and Lipschitz continuous;
\item The gradient of $\bar{f}_j^{(k)}(\bar{\bx})$ evaluated at any $\bar{\bx}_j^{(k)}$ \rev{in the interior of the Minkowski sum of $\cX_{i}^{(k)}, \,  i \in \bar{\cD}_j$}, is given by:
\begin{align} 
\nabla_{\bar{\bx}} \bar{f}_j^{(k)}|_{\bar{\bx} = \bar{\bx}_j^{(k)}} = - \bxi_j^{(k)} \label{eq:derivative_agg}
 \end{align}
where $\bxi_j^{(k)}$ is the optimal dual variable associated with constraint~\eqref{eqn:constr_sum}. \hfill $\Box$
\end{enumerate} 
\end{theorem}
  

\begin{theorem}
\label{prop:derivative_lip}
Under Assumptions~\ref{ass:agg_cost}-\ref{ass:dual_variable}, the gradient $\nabla_{\bar{\bx}} \bar{f}_j^{(k)}$ is Lipschitz continuous \rev{over the interior of the Minkowski sum of $\cX_{i}^{(k)}, \,  i \in \bar{\cD}_j$}. \hfill $\Box$
\end{theorem}


Proofs are provided in the Appendix. The results of Theorem~\ref{prop:derivative} and Theorem~\ref{prop:derivative_lip} are valid at each time instant $t_k$. These results will be utilized in the design of the real-time algorithm to update the aggregate setpoints $\bar{\bx}^{(k)}$ of groups of DERs. They will also be leveraged to establish pertinent convergence and stability claims in Section~\ref{sec:convergence}.


Let $\blambda^{(k)}, \bmu^{(k)}, \bgamma^{(k)}, \bnu^{(k)}$, and $\bzeta^{(k)}$ be the dual variables associated with constraints \eqref{mg-power-sub2}, \eqref{mg-power-sub3}, \eqref{mg-volt1}, \eqref{mg-volt2}, and \eqref{mg-current}, respectively. The Lagrangian function associated with the problem~\eqref{Pmg2} at time $t_k$ is given by:
\begin{align*} 
L^{(k)}(\bx,\bar{\bx},\bd) & := \sum_{i \in \cD} f_j^{(k)}(\bx_j) + \sum_{j \in \bar{\cD}} + \bar{f}_j^{(k)}(\bar{\bx}_j) \nonumber \\
& \hspace{-1.8cm} + \sum_{j \in \cD} \left[s^{(k)}(\blambda - \bnu)^\sfT \bM_j \bx_j  + (\bgamma - \bmu)^\sfT \bA_j \bx_j + \bzeta^\sfT \bB_j \bx_j \right] \nonumber  \\
& \hspace{-1.8cm} + \sum_{j \in \bar{\cD}} \left[s^{(k)}(\blambda - \bnu)^\sfT \bar{\bM}_j \bar{\bx}_j  + (\bgamma - \bmu)^\sfT \bar{\bA}_j \bar{\bx}_j + \bzeta^\sfT \bar{\bB}_j \bar{\bx}_j \right] \nonumber  \\
& \hspace{-1.8cm} + s^{(k)}(\blambda - \bnu)^\sfT (\mathbf{m} - \bp_{0,{\mathrm{set}}}^{(k)})      -  (\blambda + \bnu)^\sfT  E^{(k)} \mathbf{1} \nonumber \\
& \hspace{-1.8cm} + \bgamma^\sfT(\ba^{(k)} - v^{max} \mathbf{1})  + \bmu^\sfT( v^{min} \mathbf{1} - \ba^{(k)}) - \bzeta^\sfT \bi^{max} 
 \end{align*}
where $\bd := [\bgamma^\sfT, \bnu^\sfT, \blambda^\sfT, \bmu^\sfT, \bzeta^\sfT]^\sfT$ for simplicity of exposition and $\mathbf{1}$ is a vector of ones of appropriate dimensions. Consider the following regularized Lagrangian function, where  $r_p, r_d > 0$ are regularization factors: 
\vspace*{-0.1cm}
\begin{align} 
L_r^{(k)}(\bx,\bar{\bx},\bd) & := L^{(k)}(\bx,\bar{\bx},\bd) \nonumber \\
& + \frac{r_p}{2} \|\bx\|_2^2 + \frac{r_p}{2} \|\bar{\bx}\|_2^2 - \frac{r_d}{2} \|\bd\|_2^2 \label{eq:Lagrangian_aug}
 \vspace*{-0.2cm}
 \end{align}
and notice that $L_r^{(k)}(\bx,\bar{\bx},\bd) $ is strongly convex in the primal variables and strongly concave in the dual variables. Consider then the following time-varying saddle-point problem:
\begin{align} 
\max_{\bd \in \mathbb{R}_+^{2|\cM_v| +|\cM_i|  +3}} \min_{\{\bx_j\},\{\bar{\bx}_j\}} L_r^{(k)}(\bx,\bar{\bx},\bd) \label{eq:Saddle_point}
 \end{align}
and let $\bz^{(k,\star)} := [(\bx^{(k,\star)})^\sfT, (\bar{\bx}^{(k,\star)})^\sfT, (\bd^{(k,\star)})^\sfT]^\sfT$ denote the unique primal-dual optimizer of~\eqref{eq:Saddle_point}. Similar to e.g.,~\cite{opfPursuit,optimalregulationvpp2017}, the design of the online algorithm leverages appropriate modifications of online projected-gradient methods to track the time-varying solution of~\eqref{eq:Saddle_point}. Although the optimizer of \eqref{eq:Saddle_point} is expected to be different from optimizers of the original problem, 
in Section~\ref{sec:convergence} we will show that the strong convexity and concavity of $L_r^{(k)}(\bx,\bar{\bx},\bd)$ will enable the real-time algorithm to achieve Q-linear convergence. The discrepancy between $\bx^{(k,\star)}, \bar{\bx}^{(k,\star)}$ and the solution of problem $(\textrm{P1}^{(k)})$ can be bounded as shown in~\cite{Koshal11}. The point $\bz^{(k,\star)}$ is  closely related to the so-called approximate Karush-Kuhn-Tucker (KKT) conditions (see, for example,~\cite{Andreani11}).  

Let $\alpha > 0$ be a given step size. Then, given the results of Theorem~\ref{prop:derivative} and based on the regularized time-varying saddle-point formulation~\eqref{eq:Saddle_point}, the    execution
of the following steps  at each time $t_k$ defines the proposed  online algorithm. The algorithm produces power setpoints for the DERs at each time $t_k$, $k \in \mathbb{N}$.

\noindent\makebox[\linewidth]{\rule{\columnwidth}{0.4pt}}
\textbf{Real-time optimization algorithm} \\
\noindent\makebox[\linewidth]{\rule{\columnwidth}{0.4pt}}

\noindent At each $t_k$ perform the following steps:

\noindent  \textbf{[S1a]}: Collect voltage measurements $|\hat{\bv}^{(t_k)}|$ at given measurement points $\cM_v$ and perform the following updates: 
\begin{align} 
\hspace{-.3cm}\bmu^{(k+1)} & \hspace{-.1cm} = \hspace{-.05cm} \textrm{proj}_{\mathbb{R}^{|\cM_v|}_+}\hspace{-.1cm}  \left\{\bmu^{(k)} + \alpha \hspace{-.1cm} \left(v^{min} \mathbf{1} - |\hat{\bv}^{(k)}| - r_d \bmu^{(k)} \right) \right\} \label{eq:gamma} \hspace{-.2cm} \\
\hspace{-.4cm}\bgamma^{(k+1)} & \hspace{-.2cm} = \hspace{-.05cm}  \textrm{proj}_{\mathbb{R}^{|\cM_v|}_+} \hspace{-.1cm} \left\{\bgamma^{(k)} + \alpha \left(|\hat{\bv}^{(k)}| - v^{max} \mathbf{1} - r_d \bgamma^{(k)} \right) \right\} \hspace{-.2cm}
 \end{align}
\noindent  \textbf{[S1b]}: Obtain measurements or estimates of $\hat{\bi}_L^{(k)}$ on lines of interest and perform the following updates: 
\begin{align} 
\bzeta^{(k+1)} & \hspace{-.1cm} = \hspace{-.05cm}  \textrm{proj}_{\mathbb{R}^{|\cM_i|}_+} \hspace{-.05cm} \left\{\bzeta^{(k)} + \alpha \left(|\hat{\bi}_L^{(k)}| - \bi^{max} - r_d \bzeta^{(k)} \right) \right\} \hspace{-.2cm}
 \end{align}

\noindent  \textbf{[S1c]}: Collect measurements $\hat{\bp}_0^{(k)}$ at the point of common coupling  and perform the following updates: 
\begin{align} 
\blambda^{(k+1)} & \hspace{-.1cm} = \hspace{-.05cm} \textrm{proj}_{\mathbb{R}^{3}_+} \left\{\blambda^{(k)} + \alpha \hspace{-.1cm} \left(\hat{\bp}_0^{(k)} - \bp_{0,{\mathrm{set}}}^{(k)} -  E^{(k)} \mathbf{1}_3 - r_d \blambda^{(k)} \right) \right\}  \\
\bnu^{(k+1)} & \hspace{-.1cm} = \hspace{-.05cm} \textrm{proj}_{\mathbb{R}^{3}_+} \left\{\bnu^{(k)} + \alpha \hspace{-.1cm} \left(\bp_{0,{\mathrm{set}}}^{(k)} - \hat{\bp}_0^{(k)} -  E^{(k)} \mathbf{1}_3 - r_d \bnu^{(k)} \right) \right\} 
 \end{align}

\noindent  \textbf{[S2a]}: Each device $j \in \cD$ performs the following steps:

\begin{itemize}
\item[] [S2a.1] Measure output powers $\hat{\bx}_j^{(k)}$ 

\item[] [S2a.2] Update power setpoints $\bx_j^{(k+1)}$ as follows: 
\vspace*{-0.1cm}
\begin{align} 
\bx_j^{(k+1)} = & \, \textrm{proj}_{\cX^{(k)}} \left\{ \hat{\bx}_j^{(k)} - \alpha \left(\nabla_{\bx_j} f_j^{(k)}(\hat{\bx}^{(k)}_j) \right. \right. \nonumber  \\ 
& \hspace{-.35cm} +  s^{(k)}(\blambda^{(k+1)} - \bnu^{(k+1)})^\sfT \bM_j + \bzeta^{(k+1)} \bB_j \nonumber \\
& \left. \left. \hspace{-.3cm} + \,  (\bgamma^{(k+1)} - \bmu^{(k+1)})^\sfT \bA_j + r_p \hat{\bx}_j^{(k)} \right) \right\}
\end{align}

\item[] [S2a.3] If DER $j \in \cD$ has a set of discrete  setpoints, compute the implementable setpoint as:
\begin{align} 
\bepsilon_{j}^{(k)} & = \sum_{\ell = 1}^k \left(\bx_{j}^{(\ell)} - \tilde{\bx}_{j}^{(\ell)} \right) \\
\tilde{\bx}_{j}^{(k+1)} & \in \proj_{\tilde{\cX}^{(k)}_{j}} \{\bx_{j}^{(k+1)} + \bepsilon_{j}^{(k)} \} \, . \label{eq:tilde_x_der}
\end{align}

\item[] [S2a.4] Command setpoint to the DER. 

\end{itemize}

\noindent  \textbf{[S2b]}: Each DER aggregation $j \in \bar{\cD}$ performs the the following steps:

\begin{itemize}
\item[] [S2b.1] Measure aggregate output powers $\hat{\bar{\bx}}_j^{(k)}$ 

\item[] [S2b.2] Update setpoints for the aggregate powers $\bar{\bx}_j^{(k+1)}$: 
\begin{align} 
\bar{\bx}_j^{(k+1)} = & \, \textrm{proj}_{\bar{\cX}^{(k)}} \left\{ \hat{\bar{\bx}}_j^{(k)} - \alpha \left( - \bxi_j^{(k)} \right. \right. \nonumber  \\ 
& \hspace{-.35cm} +  s^{(k)}(\blambda^{(k+1)} - \bnu^{(k+1)})^\sfT \bar{\bM}_j  + \bzeta^{(k+1)}  \bar{\bB_j}  \nonumber \\
& \left. \left. \hspace{-.3cm} + \,  (\bgamma^{(k+1)} - \bmu^{(k+1)})^\sfT \bar{\bA}_j + r_p \hat{\bar{\bx}}_j^{(k)} \right) \right\}
\end{align}

\item[] [S2b.3] Given the aggregate powers $\bar{\bx}_j^{(k+1)}$, compute the setpoints $\{\bx_i \in \cX_i^{(k)}\}_{i \in \bar{\cD}_j}$ of the individual DERs $\bar{\cD}_j$ and the new vector $\bxi_j^{(k+1)}$ by solving the saddle-point problem:
\begin{equation}
\max_{\bxi} \min_{\{\bx_i \in \cX_i^{(k)}\}_{i \in \bar{\cD}_j}} \sum_{i \in \bar{\cD}_j} f_i^{(k)}(\bx_i) + \bxi^\sfT \Big( \sum_{i \in \bar{\cD}_j} \bx_{i} - \bar{\bx}_j^{(k+1)}\Big).
\end{equation}

\item[] [S2b.4] If DER $j \in \bar{\cD}_j$ has a set of discrete  setpoints, compute the implementable setpoint as:
\begin{align} 
\bepsilon_{j}^{(k)} & = \sum_{\ell = 1}^k \left(\bx_{j}^{(\ell)} - \tilde{\bx}_{j}^{(\ell)} \right) \\
\tilde{\bx}_{j}^{(k+1)} & \in \textrm{proj}_{ \tilde{\cX}^{(k)}_{j}} \left\{\bx_{j}^{(k+1)} + \bepsilon_{j}^{(k)}  \right\} \, . \label{eq:tilde_x}
\end{align}

\item[] [S2b.5] Command setpoints to the DERs. 

\end{itemize}

\noindent  \textbf{[S3]}: Go to \textbf{[S1]}.

\vspace{.2cm}

\rev{The following remarks are in order:}

\begin{figure}
  \centering
  \includegraphics[width=1.0\columnwidth]{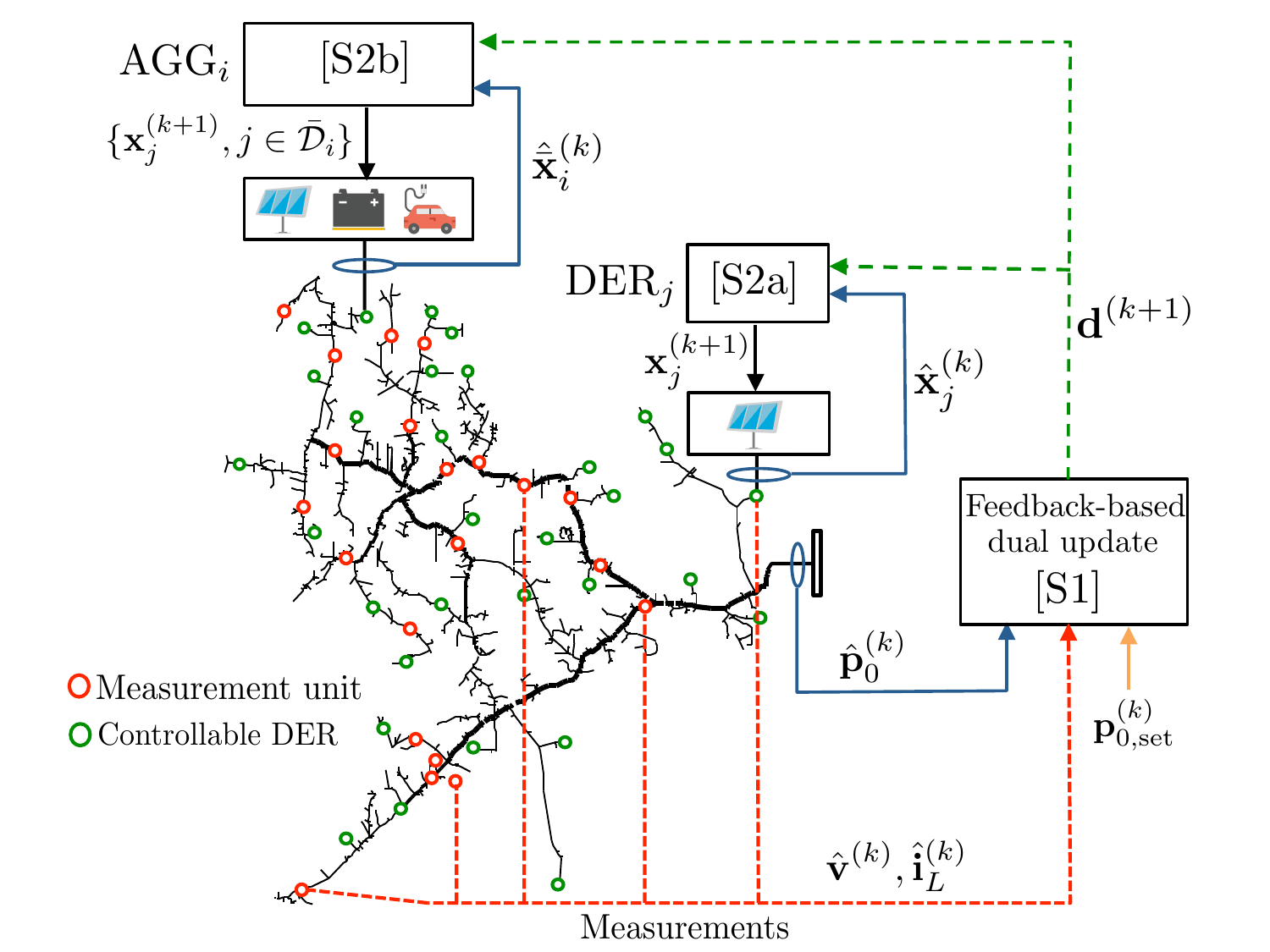}
\caption{Feedback-based online algorithm: distributed implementation. } \label{fig:F_algorithm}
\vspace{-.5cm}
\end{figure}

\begin{enumerate}
    \item Notice that the feedback is utilized in the algorithm in steps~\textbf{[S1]} and ~\textbf{[S2]} in the form of measurements of voltages, currents, and power flows at the point of common coupling; these measurements replace the corresponding analytical expressions. By virtue of this approach, challenges (c2)-(c3) are resolved and, in particular, measurement of the state of uncontrollable devices is not required. \rev{We note that in its straightforward implementation, the algorithm requires measurements at nodes where there a corresponding constraint is imposed; however, if real-time state estimation procedure is available, the measurements can be replaced with the estimated state.} 
    
    \item The real-time algorithm affords a distributed implementation as shown in  Fig.~\ref{fig:F_algorithm}. Once measurements  $\hat{\bv}^{(k)}$, $\hat{\bi}_L^{(k)}$, and $\hat{\bp}_0^{(k)}$ are acquired, step~\textbf{[S1]} is performed at the utility/aggregator, which subsequently broadcasts the dual variables $\bd^{(k+1)}$. Steps~\textbf{[S2a]} and~\textbf{[S2b]} are implemented locally at individual DERs and aggregations of DERs (in the Fig.~\ref{fig:F_algorithm}, AGG stands for aggregation), respectively. 
    
    \item \rev{Note that the sub-steps \textbf{[S1a], [S1b]}, and \textbf{[S1c]} can be carried out in parallel at the utility/aggregator.}

    \item \rev{ The framework is flexible enough so that the entire computation can be performed centrally. However, there are obviously disadvantages to this. First, this  creates a single point failure: if the central entity fails, the entire scheme stops working; however, in the distributed implementation, the local controllers can still compute (sub-)optimal setpoints using perhaps outdated network-wide information (i.e., outdated Lagrange multipliers). Second, the centralized implementation requires point-to-point communication with every DER to communicate individual setpoints; on the other hand, distributed implementation advocated here only requires broadcast communication -- every DER receives the same Lagrange multiplier. This is especially important when the number of measurement points (or, the points where we want to impose constraints) is much smaller than the number of DERs in the system. Finally, the privacy argument applies here: if we implement the algorithm centrally, the DERs have to reveal their private information (e.g., preferences and feasible regions).}
    
    \item \rev{The steps \eqref{eq:tilde_x_der} and \eqref{eq:tilde_x} represent the implementation based on the error-diffusion algorithm. In particular, the accumulate error between the continuous setpoint $\bx_j^{(k)}$ and the discrete implementation $\tilde{\bx}_j^{(k)}$ is computed and used in the modified projection steps \eqref{eq:tilde_x_der} and \eqref{eq:tilde_x} to obtain the next implementable (discrete)  setpoint; see \cite{errDiff} for further details.
    }
    Finally, notice that steps~\eqref{eq:tilde_x_der} and~\eqref{eq:tilde_x} involve the solution of a localized nonconvex program to compute implementable commands.

\end{enumerate}

The ability of the algorithm to track the optimizers $\bz^{(k, \star)}$ of~\eqref{eq:Saddle_point} is analytically established next.

\section{Performance Analysis}
\label{sec:convergence}

\rev{We next analyze the proposed algorithm under the assumption of \emph{synchronous} updates in steps \textbf{[S1]} and \textbf{[S2]} above. The analysis of the asynchronous case can be carried out similarly on expense of heavier notation and further assumptions; see for example \cite{Bernstein2018}.}

We start by stating the following assumption regarding problem~\eqref{Pmg2}.  

\begin{assumption} \label{asm:feasibility}
Problem~\eqref{Pmg2} is feasible and Slater's condition holds at each time $t_k$, $k \in \mathbb{N}$. 
\hfill $\Box$
\end{assumption}

Assumption~\ref{asm:feasibility} implies  that there exists a power flow solution that adheres to voltage and ampacity limits. When the distribution network is required to follow a setpoint at the point of common coupling,  Assumption~\ref{asm:feasibility} presumes that the setpoint is feasible. Feasibility of the power flow solutions (with and without setpoints for the active and reactive power at the substation) can be assessed by solving suitable optimization problems at a slower time scale; see, for example, the optimization approaches proposed in~\cite{Nosair15,Meliopoulos16}. \\

Regarding the temporal variability of problem~\eqref{Pmg2}, we introduce the following quantity to capture the variation of the optimal solution trajectory over time: 
\rev{
\begin{align} 
\label{eqn:sigma}
\sigma^{(k)} &:= \|\bz^{(k+1,\star)} - \bz^{(k, \star)}\|_2, \quad \sigma := \sup_{k \geq 1} \sigma^{(k)}.
\end{align} 
}
For sufficiently small sampling intervals $h$, $\sigma$ can be interpreted as a bound on the norm of the gradient of the optimal solution trajectory $\{\bz^{(k,\star)}\}_{k \in \mathbb{N}}$ with respect to time. In the context of~\eqref{Pmg2}, $\sigma$ depends on the variability of the cost function, non-controllable loads, as well as available powers from the renewable-based DERs.  

Next, since models~\eqref{eq:approx_v}--\eqref{eq:approx_p} are linear and the sets $\{\cX_j^{(k)}\}$ and $\{\bar{\cX}_j^{(k)}\}$ are compact (cf.~Assumption \ref{ass:agg_cost}), there exist constants $G_v < + \infty$, $G_0 < + \infty$, and $G_L < + \infty$ such that, for every $k \in   \mathbb{N}$,
$$\|\nabla_{[\bx, \bar{\bx}]} |\tilde{\bv}^{(k)}(\bx, \bar{\bx})| \|_2 \leq G_v ,
\|\nabla_{[\bx, \bar{\bx}]} \tilde{\bp}^{(k)}_0(\bx, \bar{\bx}) \|_2 \leq G_0 ,$$ 
$$\|\nabla_{[\bx, \bar{\bx}]} |\tilde{\bi}^{(k)}_L(\bx, \bar{\bx})| \|_2 \leq G_L .$$
For future developments, define $G := \max\{G_v, G_0, G_L\}$. Further, notice that from Assumption~\ref{ass:agg_cost} and Theorem~\ref{prop:derivative_lip}, the gradient map $\bg^{(k)}(\bx, \bar{\bx}) := [\nabla_{\bx_1}^\sfT f_1^{(k)}(\bx_1), \ldots, \nabla_{\bx_{|\cD|}}^\sfT f_{|\cD|}^{(k)}(\bx_{|\cD|}),$  $\nabla_{\bar{\bx}_1}^\sfT \bar{f}_1^{(k)}(\bar{\bx}_1), \ldots, \nabla_{\bar{\bx}_{|\bar{\cD}|}}^\sfT \bar{f}_{|\bar{\cD}|}^{(k)}(\bar{\bx}_{|\bar{\cD}|})]^\sfT$ is 
 Lipschitz continuous with a given constant $L^{(k)}$ over the set $\cX^{(k)} := \cX_1^{(k)} \times \ldots \times \cX_{|\cD|}^{(k)} \times \bar{\cX}_1^{(k)} \times \ldots \times \bar{\cX}_{|\bar{\cD}|}^{(k)}$. Let $L := \sup \{L^{(k)}\}$, so that 
\begin{align}
\|\bg^{(k)}(\bx, \bar{\bx}) - \bg^{(k)}(\bx^\prime, \bar{\bx}^\prime)\|_2 \leq L \|\bx - \bx^\prime\|_2 
\end{align}
 for all $\bx, \bx^\prime \in \cX^{(k)}$, $\bar{\bx}, \bar{\bx}^\prime \in \bar{\cX}^{(k)}$, and $t_k$, $k \in \mathbb{N}$.

Define the errors introduced by the measurement noise and modeling mismatches (i.e., discrepancy between the nonlinear AC power-flow equations and the linearized model, as well as possible inaccurate knowledge of the admittance matrix) as follows: 
\begin{align*} 
e^{(k)}_x & := \left\| 
\left[
\begin{array}{c}
\bx^{(k)} \\
\bar{\bx}^{(k)}
\end{array}
\right]
- 
\left[
\begin{array}{c}
\hat{\bx}^{(k)} \\
\hat{\bar{\bx}}^{(k)}
\end{array}
\right]
\right\|_2  \\ 
e^{(k)}_0 & := \|\tilde{\bp}_0(\bx^{(k)}, \bar{\bx}^{(k)})  - \hat{\bp}_0^{(k)}\|_2 \\
e^{(k)}_v & := \| |\tilde{\bv}^{(k)}(\bx^{(k)}, \bar{\bx}^{(k)}) - |\hat{\bv}^{(k)}| \|_2  \\
e^{(k)}_L & := \| |\tilde{\bi}_L^{(k)}(\bx^{(k)}, \bar{\bx}^{(k)}) - |\hat{\bi}_L^{(k)}| \|_2 
\end{align*}
where we recall that $\hat{\bv}^{(k)}$, $\hat{\bi}_L^{(k)}$, and $\hat{\bp}_0^{(k)}$ are measurements (or pseudo-measurements). The following assumption is made.

\begin{assumption} \label{asm:bounded}
There exist finite constants $e_x$, $e_0$, $e_v$, and $e_L$ such that $e_x^{(k)} \leq e_x$, $e_0^{(k)} \leq e_0$, $e_v^{(k)} \leq e_v$, and $e_L^{(k)} \leq e_L$  for all $t_k$; that is, the errors are uniformly bounded in time.  \hfill $\Box$
\end{assumption}

\rev{As discussed in Section~\ref{sec:modelDER}, DERs are presumed to be equipped with embedded  controllers that drive the output powers to the commanded setpoints. }
If the time constant of the controllers is longer than $h$, Assumption~\eqref{asm:bounded} bounds the discrepancy between the sampled output power and the commanded setpoint. For future developments,  define the vector $\be^{(k)} := [(L+r_p) e_x^{(k)},  \mathbf{1}_2^\sfT e_v^{(k)}, \mathbf{1}_2^\sfT  e_0^{(k)}, e_L^{(k)} ]^\sfT$, and notice from Assumption~\ref{asm:bounded} that $\|\be^{(k)}\|_2 \leq e$, $e := \sqrt{(L+r_p)^2 e_x^2 + 2 e_v^2 + 2 e_0^2 + e_L^2}$.  

Let $\bz^{(k)} := [(\bx^{(k)})^\sfT, (\bar{\bx}^{(k)})^\sfT, (\bd^{(k)})^\sfT]^\sfT$ collect the primal and dual variables produced by the real-time algorithm at time $t_k$. Based on Assumptions~\ref{ass:agg_cost}--\ref{asm:bounded}, the main convergence results are established next.

\begin{theorem}
\label{theorem.inexact}
Consider the sequence $\{\bz^{(k)}\}$ generated by the algorithm~\eqref{eq:gamma}--\eqref{eq:tilde_x}. The distance between $\bz^{(k)}$ and the primal-dual optimizer $\bz^{(k,\star)}$ at time $t_k$ can be bounded as:
\begin{align} 
\|\bz^{(k)} - \bz^{(k, \star)}\|_2 \leq & c(\alpha)^k \|\bz^{(0)} - \bz^{(0, \star)}\|_2 \nonumber \\ 
& \hspace{-2.5cm} + \sum_{\ell = 0}^{k-1} c(\alpha)^\ell \left(e^{(k-\ell-1)}_x + \alpha \|\be^{(k-\ell-1)}\|_2 + \sigma^{(k-\ell-1)} \right)
\label{eqn:inst_bound}
\end{align}
where 
\begin{align}
c(\alpha) := & [1 - 2 \alpha \min\{r_p,r_d\} + \alpha^2(L+r_p+5G)^2  \nonumber \\
& + 5 \alpha^2 (G + r_d)^2]^\frac{1}{2} \, \label{eq:kappa}
\end{align}
and $\sigma^{(k)}$ is defined in~\eqref{eqn:sigma}. 
\hfill $\Box$
\end{theorem}

\begin{corollary} \label{cor:asym}
If $c(\alpha) < 1$, then the sequence  $\{\bz^{(k)}\}$ converges Q-linearly to $\{\bz^{(k, \star)}\}$ up to an asymptotic error bound given by:
\begin{align} 
\limsup_{k\to\infty} \|\bz^{(k)} - \bz^{(k, \star)}\|_2 &\leq \frac{\Delta}{1 - c(\alpha)}, \label{eqn:asym_bound_apr}
\end{align}
where $\Delta := e_x + \alpha e + \sigma$. \hfill $\Box$
\end{corollary}

\vspace{.1cm}

Notice first that the condition $c(\alpha, r_p, r_d) < 1$ is satisfied if:
\begin{align} 
\alpha < \frac{\min\{r_p, r_d\}}{(L+r_p+5G)^2 +  5(G + r_d)^2}.
\end{align} 
\rev{Also, observe that the value of $\Delta$ (and hence knowledge of $e_x, e$, and $\sigma$) is not required in order to satisfy the condition of Corollary \ref{cor:asym}.}

The bound~\eqref{eqn:inst_bound} provides a characterization of the discrepancy between $\bz^{(k, \star)}$ and $\bz^{(k)}$ at each time $t_k$. On the other hand,  the asymptotic bound~\eqref{eqn:asym_bound_apr} depends on the underlying dynamics of the distribution system through $\sigma$ and on the measurement errors through $e$. The result~\eqref{eqn:asym_bound_apr}  can also be interpreted as input-to-state stability, where the optimal trajectory $\{\bz^{(k, \star)}\}$ of the time-varying problem~\eqref{Pmg2} is taken as a reference. When $e = 0$ and $\sigma = 0$, the algorithm converges to the solution of the static optimization problem~\eqref{eq:Saddle_point}. The proof of Theorem~\ref{theorem.inexact} is sketched in Appendix~\ref{sec:proofthm_inexact}.

We conclude the section by stating \rev{for completeness} a result from \cite{errDiff} establishing  average tracking properties for the updates~\eqref{eq:tilde_x_der} and~\eqref{eq:tilde_x}. 

\rev{
\begin{theorem}[Theorem 2 in \cite{errDiff}] \label{thm:err_diff}
   For each DER $j$ with nonconvex operational region $\cX_j^{(k)}$ there exists a finite constant $E_{j}$ such that $\|\bepsilon_{j}^{(k)}\|_2 \leq E_{j}$ for all $k$. Consequently, 
\begin{align}
\label{eq:error_discrete}
   \left \| \frac{1}{k}\sum_{\ell = 1}^k \bx_{j}^{(\ell)} - \frac{1}{k}\sum_{\ell = 1}^k \tilde{\bx}_{j}^{(\ell)} \right \|_2 \leq \frac{E_{j}}{k}
\end{align}
and $\|\bx_{j}^{(k)} - \tilde{\bx}_{j}^{(k)}\|_2 \leq 2 E_{j}$ for all $k$. \hfill $\Box$
   \end{theorem}
}

Representative numerical experiments using real data are presented in the next section.

\section{Experiments on a Real System}
\label{sec:results}

\begin{figure*}[t]
  \centering
  \includegraphics[width=1.5\columnwidth]{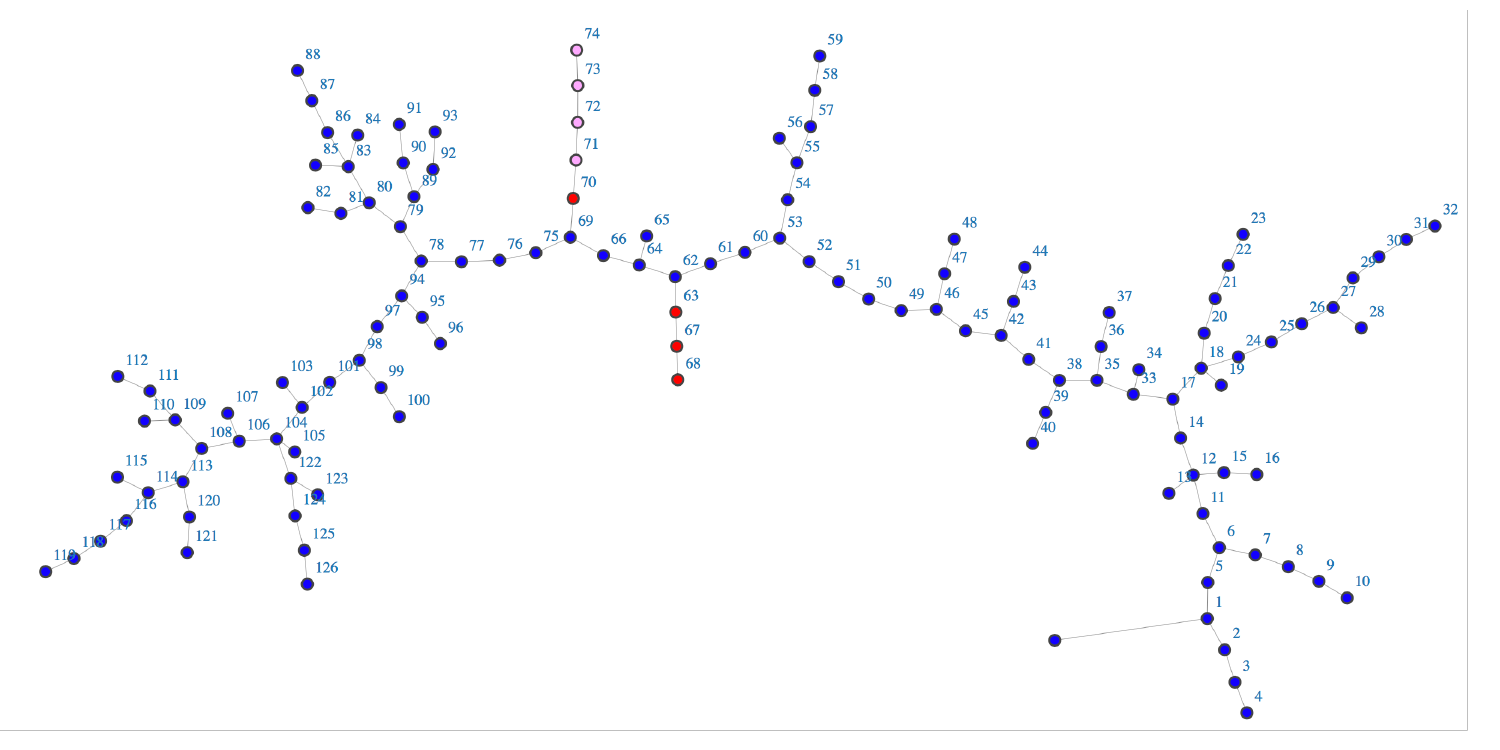}
\vspace{-.2cm}
\caption{Anonymized diagram of the SCE feeder utilized in the numerical experiments. Nodes are color-coded based on the number of phases: blue nodes have three phases, red nodes have two phases, and pink nodes are single-phase.} \label{fig:F_circuit}
\vspace{-.2cm}
\end{figure*}

The proposed real-time algorithm is tested using data of a real distribution feeder located within the territory of Southern California Edison (SCE)
As shown in the anonymized diagram in Fig.~\ref{fig:F_circuit}, 
this distribution feeder features 126 multiphase nodes (excluding the substation), with a total of 366 single-phase points of connection. Wye and delta connections are present at different nodes of the feeder. The feeder has a nominal line-line voltage of 12kV, and it has three phases in all the nodes except the following ones: nodes 63, 67, 68, and 70 have only phases b and c, and nodes 71, 72, 73, and 74 have only phase c (see node numbering in Appendix \ref{sec:data}). The feeder serves 362 customers, with a mix of residential,  commercial, and industrial facilities.  In the numerical experiments, controllable assets include photovoltaic (PV) systems, energy storage systems, and electric vehicles (EVs). Load and irradiance data have a granularity of 6 seconds; to achieve a granularity of 1 second, the time series were interpolated. It follows that the target optimization problem~\eqref{Pmg2}--\eqref{eq:aggregation} changes every second. 

\begin{figure}[t]
  \centering
  \includegraphics[width=1.0\columnwidth]{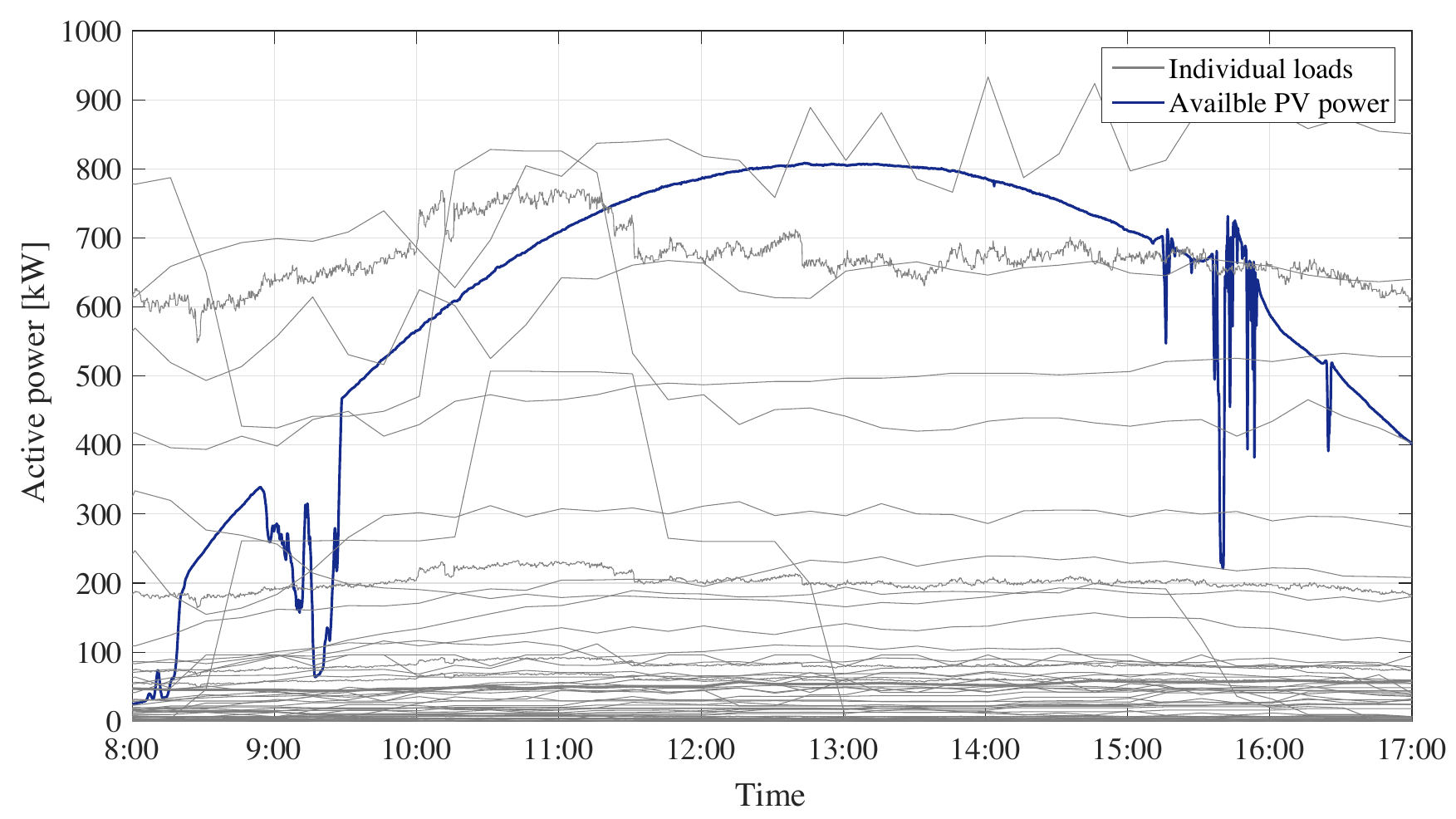}
\vspace{-.8cm}
\caption{Trajectories of individual non-controllable loads and power available from PV systems with capacity of 1MW.} \label{fig:F_loads}
\vspace{-.3cm}
\end{figure}
\begin{figure}[t]
  \centering
  \includegraphics[width=1.0\columnwidth]{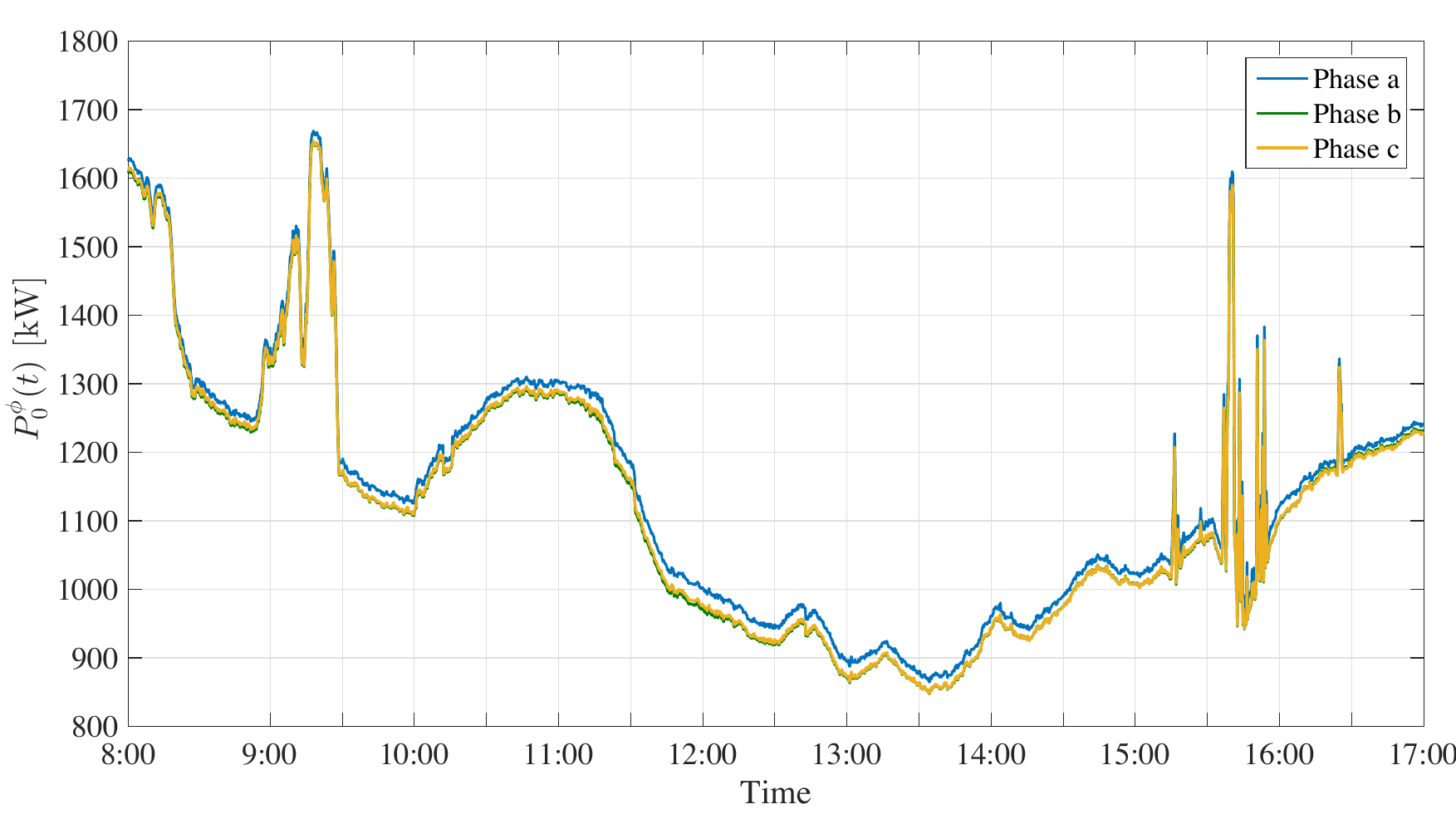}
\vspace{-.8cm}
\caption{Trajectories of the power at the substation when devices are controlled as business-as-usual.} \label{fig:F_net_load}
\vspace{-.3cm}
\end{figure}

\rev{It is worth pointing out that this feeder is ``stiff'' and includes only the modelling of the primary side of distribution transformers. Therefore, voltage violations in this feeder are less visible. The main reason for choosing this test case is that this is a real feeder from SCE with real data. The downside is that SCE does not provide any models on the secondary side, where the possibility of the violation is more likely.}

The algorithm was coded in Matlab. For a given set of net power injections/consumptions at each of the nodes of the feeder, a power flow solution was obtained using OpenDSS. 

The location of the PV and battery systems in the feeder, along with their capacities and connections types, can be found in Appendix \ref{sec:data}. In case of multiple PV systems at a node, the devices are aggregated and jointly controlled. Three-phase systems are presupposed to operate in a balanced mode. The operating region of the inverters that accompany PV systems and energy storage systems is in the form of~\eqref{eq:set}.  On the other hand, level-2 charging stations for EVs are presupposed, with discrete charging levels of 10, 20, 40, 60, 80, and 100\% of the maximum charging capability of 7.2 kW. EVs are located as follows: 5 EVs at node 9, 2 at node 29, and 3 at node 90. The batteries of the EVs have sized of 60, 80, and 130 kWh, and a minimum charging rate is set for the EVs so that they can be fully charged at the time specified by the drivers.

Table~\ref{tab:pv_systems} in Appendix~\ref{sec:data} lists the locations of the PV systems, along with their capacities and connections types; for the latter, the symbol $xY$ refers to a wye connection on $x = |\cP_j|$ phases and $x\Delta$ refers to  $x$ delta connections between the available phases. In case of multiple PV systems at a node, the devices are aggregated and jointly controlled. Three-phase systems are presupposed to operate in a balanced mode. Table~\ref{tab:battery_systems} lists the locations of the battery systems, along with the capacity of the inverters, the maximum state of charge, and the connection type. Similarly to PV systems, three-phase batteries operate in a balanced mode. The operating region of the inverters that accompany PV systems and energy storage systems is in the form of~\eqref{eq:set}.  On the other hand, level-2 charging stations for EVs are presupposed, with discrete charging levels of 10, 20, 40, 60, 80, and 100\% of the maximum charging capability of 7.2 kW. EVs are located as follows: 5 EVs at node 9, 2 at node 29, and 3 at node 90. The batteries of the EVs have sized of 60, 80, and 130 kWh, and a minimum charging rate is set for the EVs so that they can be fully charged at the time specified by the drivers. 

The trajectories of individual non-controllable loads and the power available from a PV system with an inverter capacity of 1MW  for a selected day in September 2016 are illustrated in Fig.~\ref{fig:F_loads}. The power available from other PV systems is a scaled version of the trajectory shown in Fig.~\ref{fig:F_loads}. It can be seen that the selected day is mostly sunny with clear sky; however, clouds introduced a significant variability in the available power from 9:00 to 9:30 and in the afternoon from 15:15 to 16:00. 
Fig.~\ref{fig:F_net_load} illustrates the trajectories of the active power at the substation when PV systems are operated at the maximum available power and batteries are not utilized. It can be seen that the feeder is unbalanced, in the sense that there is a discrepancy between the  power of phase $a$ and that of the remaining phases at the substation. 
We note that the majority of the controllable assets are three-phase (balanced) with delta connections; hence, it is not possible to balance the operation of the feeder (that is, ensure that the net powers at the three phases of the substation are equal at each point in time). 

In the first test, the algorithm is evaluated during the sunny period of the day; in a second test, we test the algorithm during cloudy periods to assess whether it can cope with uncertain (and fast-changing) weather conditions. The PV-related cost functions are set to $(P_{\textrm{av},i}^{(k)} - P_{j}^{(k)})^2 +  (Q_{j}^{(k)})^2$ for three-phase PV systems (with $P_{\textrm{av},j}^{(k)}$ denoting the maximum real power available) and $100(P_{\textrm{av},i}^{(k)} - P_{j}^{(k)})^2 +  10 (Q_{j}^{(k)})^2$ for smaller-size single-phase PV systems. For batteries, $f_i^{(k)}$ and $\bar{f}_i^{(k)}$ are set to $(P_{j,\phi}^{(k)})^2 +  (Q_{j,\phi}^{(k)})^2$, and for the EVs we have that $100 (P_{j}^{(k)} - P_{\textrm{max},i})^2$ where $P_{\textrm{max},i}$ is the maximum charging rate. With this setting, the DER will be incentivized to provide services to the grid, while minimizing the power curtailed from the PV systems and  the deviation from a predetermined (dis)charging profiles for the batteries. The  stepsize is set as $\alpha = 0.2$ and the regularization parameters are $r_p = 10^{-3}$ and $r_d = 10^{-4}$. One step of the algorithm is run every $1$ second. PV and battery inverters are presumed to follow a first-order response with a time constant of $0.25$ seconds (hence, a settling time of approximately $1.25$ seconds). Communication delays are set to $0.1$ seconds per link. Voltage limits are set to $0.95$ pu and $1.05$ pu. 

\begin{figure}[t]
  \centering
  \includegraphics[width=1.0\columnwidth]{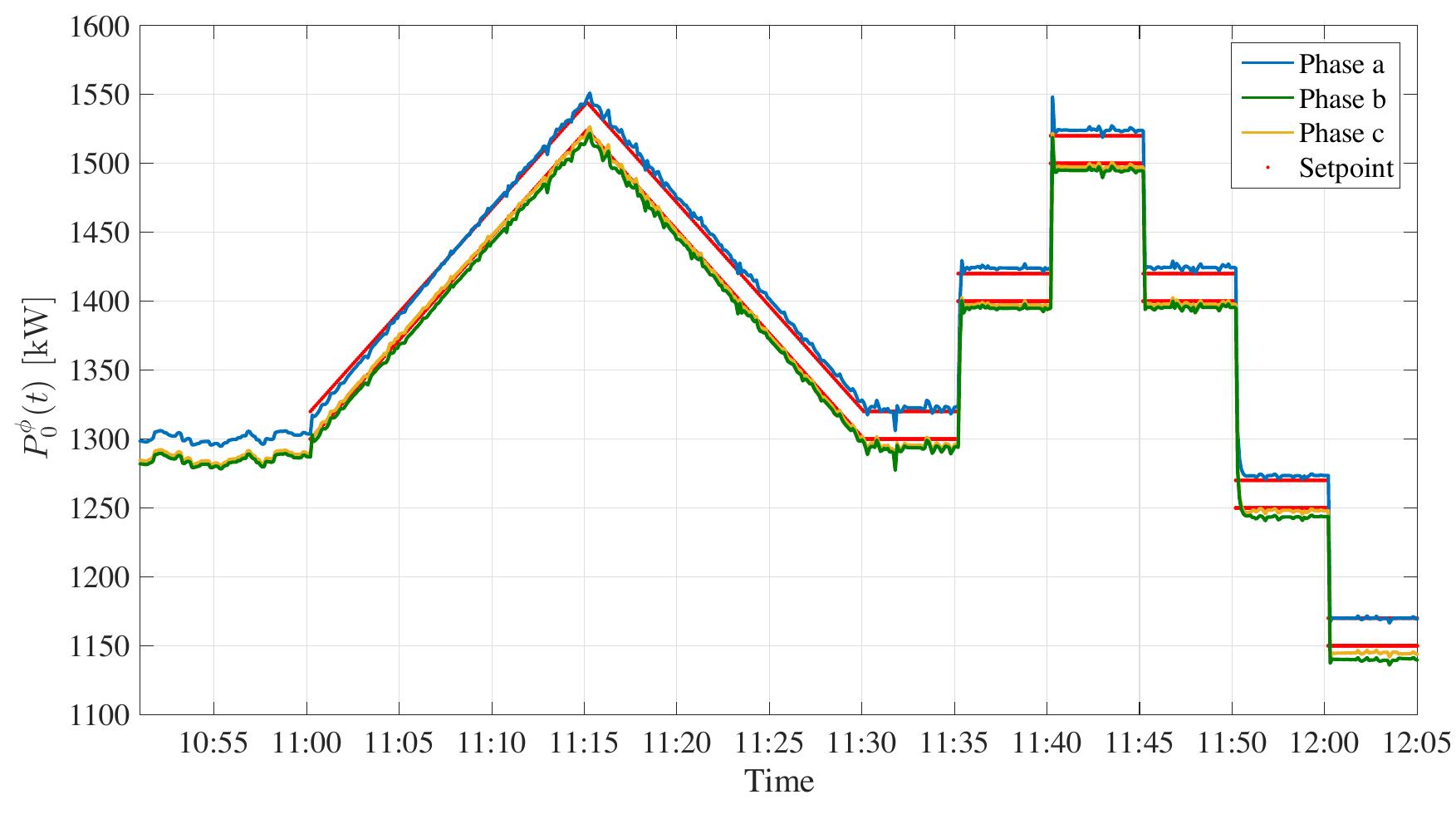}
\vspace{-.8cm}
\caption{Tracking of setpoints for the power at the substation, while respecting voltage limits. Clear sky case. } \label{fig:F_track_1}
\vspace{-.2cm}
\end{figure}
\begin{figure}[t]
  \centering
  \includegraphics[width=1.0\columnwidth]{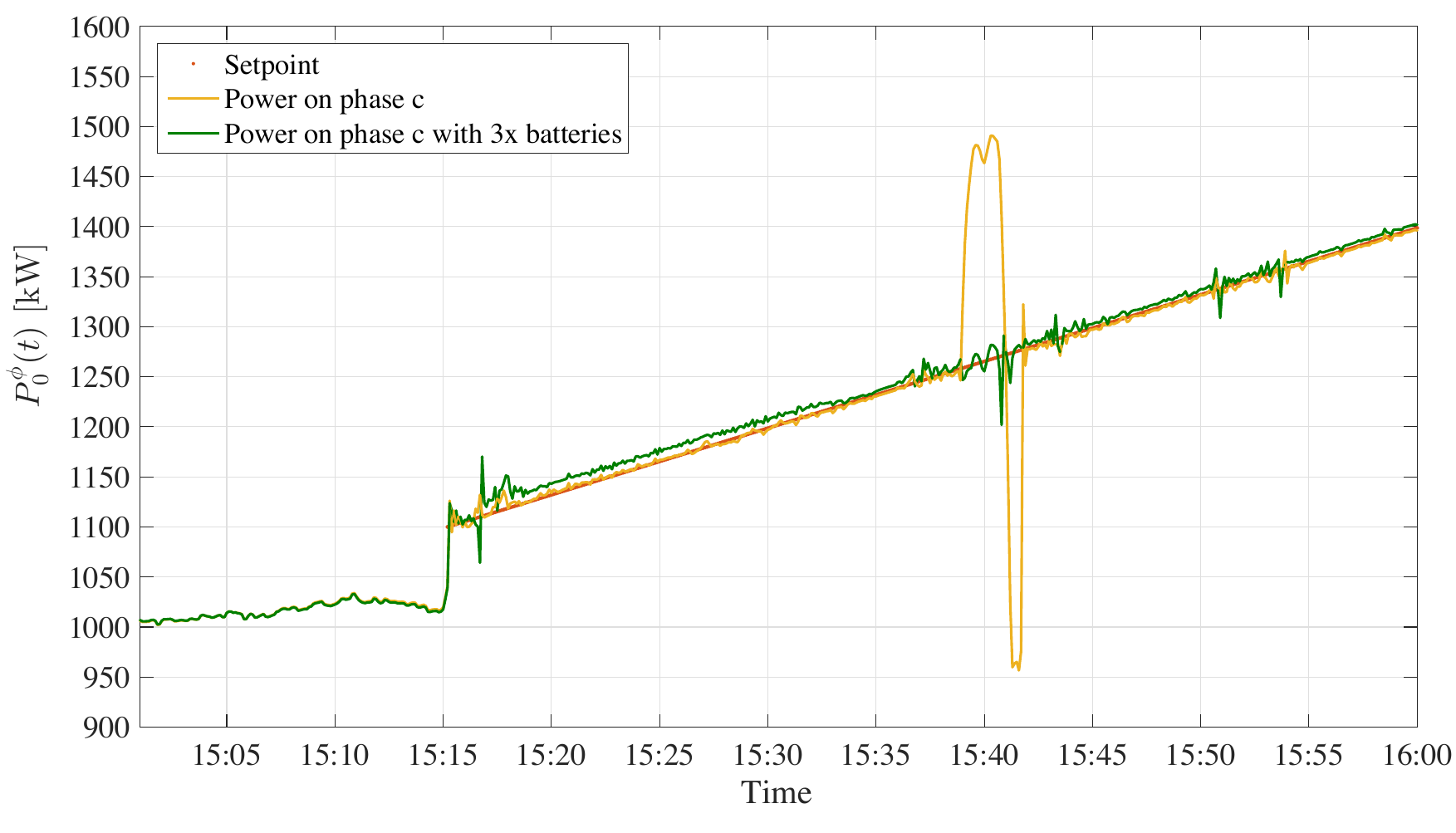}
\vspace{-.8cm}
\caption{Tracking of setpoints for the power at the substation, while respecting voltage limits. Cloudy sky case.} \label{fig:F_track_2}
\vspace{-.4cm}
\end{figure}
\begin{figure}[t]
  \centering
  \vspace{-.2cm}
  \includegraphics[width=1.0\columnwidth]{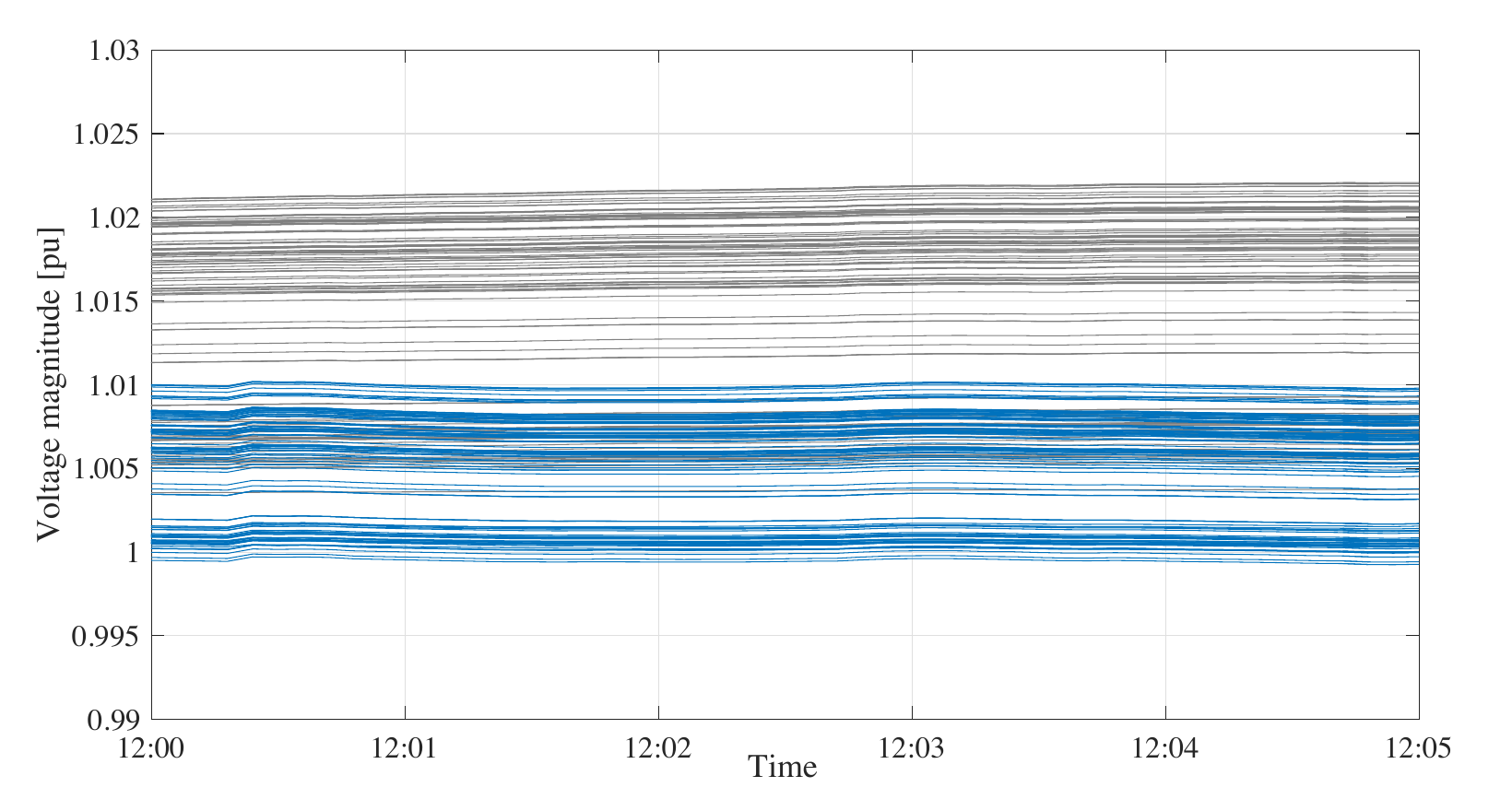}
\vspace{-.8cm}
\caption{Voltages magnitude, for a representative time slot. Grey: upper limit set to $1.05$ pu. Blue: $1.01$ pu. } 
\label{fig:F_voltage}
\vspace{-.6cm}
\end{figure}

In the test cases, we control the DERs in order to track a given trajectory of setpoints $\bp_{0,{\mathrm{set}}}^{(k)}$ at the substation, while ensuring that voltages are within limits.
Fig.~\ref{fig:F_track_1} illustrates the tracking performance of the real-time algorithm from 11:00 to 12:00, where the majority of the problem variability is introduced by non-controllable devices. The red trajectory corresponds to setpoints $\bp_{0,{\mathrm{set}}}^{(k)}$ (which are different across phases to acknowledge the unbalance operation of the feeder), while the powers on the three phases of the substation are color-coded in blue (phase a), green (phase b) and orange (phase c). It can be seen that with the proposed algorithm the power at the substation closely track ramping signals~\cite{Caiso12} as well as step changes in the setpoints. Similar tracking results are shown in Fig.~\ref{fig:F_track_2}, where we considered the time interval from 15:30 to 16:00, where the the overall power available from the PV systems is varying very fast. In the uncontrolled case in Fig.~\ref{fig:F_net_load}, the variation in the power available from the PV systems translated into spikes in the power at the substation of magnitude up to $1.8$ MW (summed across phases). On the other hand, the algorithm is capable of leveraging energy storage system to lower the power swing. In addition, Fig.~\ref{fig:F_track_2} reports the tracking result in the case where we increase the capacity of the batteries of 3x. In this case, the algorithm is capable of completely copying with the PV variability (only the phase c is reported to facilitate the comparison between the two cases).  

In the previous experiments, voltage magnitudes were well within limits because the feeder is stiff. To test the ability of our method to regulate voltages, we increase the capacity of the PV systems of 3x to create reverse power flow conditions, and we lower the upper limit from $1.05$ to $1.01$ pu. Fig.~\ref{fig:F_voltage} illustrates the ``cloud'' of voltages magnitudes across the system for a representative time slot; the blue trajectories represent the voltages magnitudes when the upper limit is $1.01$ pu. It can be seen that the algorithm is capable of regulating voltages while driving the power at the substation to specific setpoints. 

\vspace{-0.3cm}

\section{Concluding Remarks}
\label{sec:conclusions}

This paper developed a distributed algorithm  for real-time optimization of DERs. The proposed framework optimizes the operation of both DERs that are individually controllable and groups of DERs  at an electrical point of connection that are jointly controlled, and it enables (groups of) DERs to pursue given performance objectives while adjusting their (aggregate) powers to respond to services requested by grid operators and to maintain electrical quantities within engineering limits. The design of the algorithm leveraged a  time-varying bi-level problem formulation capturing various performance objectives and engineering constraints, and a feedback-based online implementation of primal-dual projected-gradient methods. The resultant feedback-based online algorithm can cope with inaccuracies in the distribution-system modeling, it avoids pervasive metering to gather the state of non-controllable resources, and it naturally lends itself to a distributed implementation. Analytical stability and convergence claims were established in terms of tracking of the solution of the formulated time-varying  optimization problem. Future efforts will look at extending the technical findings to time-varying nonconvex problems.



\vspace{-.5cm}

\appendix

\subsection{Proof of Proposition \ref{thm:sum1}}
\label{sec:proofthm1}
The proof follows from the fact that the Minkowski sum between a circle $\cC := \{(x, y): \, x^2 + y^2 \leq r^2\}$ and an interval $\cI := \{(x, y): \, a \leq x \leq b, y = 0\}$ is given by:
\[
\cC + \cI = \{(x, y): \, a - r \leq x \leq b + r, -h(x) \leq y \leq h(x) \},
\]
where $h(x)$ is a concave function given by
\[
h(x) :=
\begin{cases}
r, & x \in [a, b] \\
\sqrt{r^2 - (x - a)^2}, & x \in [a-r, a) \\
\sqrt{r^2 - (x - b)^2}, & x \in (b, b + r] . 
\end{cases}
\]

\vspace{-0.3cm}

\subsection{Proof of Proposition \ref{thm:sum2}}
\label{sec:proofthm2}

We first prove the outer approximation~\eqref{eq:sum_outer}. Let $\bar{\bx} = [P,Q]^\sfT \in \cX(\underline{p}_1, \overline{p}_1, r_1) \oplus \cX(\underline{p}_2, \overline{p}_2, r_2)$; that is, $\bar{\bx} = \bx_1 + \bx_2$ for two DERs $\bx_1  = [P_1,Q_1]^\sfT \in \cX(\underline{p}_1, \overline{p}_1, r_1)$ and $\bx_2 = [P_2,Q_2]^\sfT \in \cX(\underline{p}_2, \overline{p}_2, r_2)$. It can be readily shown that $ \underline{p}_1 + \underline{p}_2 \leq P = P_1 + P_2 \leq \overline{p}_1 + \overline{p}_2$; further,  we  have that:
\begin{align*}
&P^2 + Q^2 = (P_1+P_2)^2 + (Q_1 + Q_2)^2 \\
&=P_1^2 + Q_1^2 + P_2^2 + Q_2^2 + 2\bx_1^\sfT \bx_2 \leq r_1^2 + r_2^2 + 2 \|\bx_1\|_2 \|\bx_2\|_2\\
&\leq r_1^2 + r_2^2 + 2 r_1 r_2 = (r_1 + r_2)^2
\end{align*}
where the first inequality follows by the fact that $\bx_i \in \cX(\underline{p}_i, \overline{p}_i, r_i)$, $i = 1, 2$, and from the Cauchy-Schwarz inequality; the second inequality holds again by construction. It follows that $\bar{\bx} \in \cX(\underline{p}_1 + \underline{p}_2, \overline{p}_1 + \overline{p}_2, r_1 + r_2)$ whenever $\bx_1 \in \cX(\underline{p}_1, \overline{p}_1, r_1)$ and $\bx_2 \in \cX(\underline{p}_2, \overline{p}_2, r_2)$; the converse is, however, not necessarily true. 

We next prove the inner  approximation~\eqref{eq:sum_inner}. Let $\bx = [P, Q]^\sfT \in \cX(\underline{p}_1 + \underline{p}_2, \overline{p}_1 + \overline{p}_2, \rho)$. Clearly,
\begin{align}
&\exists P_1, P_2: \, \underline{p}_1 \leq P_1 \leq \overline{p}_1, \, \underline{p}_2 \leq P_2 \leq \overline{p}_2, \, P_1 + P_2 = P.
\end{align}
Using the definitions of $\alpha$ and $\beta_i$, one can verify that
\begin{align} 
P_i^2 \leq \beta_i, \, i = 1, 2, \label{eqn:betaBound}\\
(P_1 + P_2)^2 \geq \alpha.\label{eqn:alphaBound}
\end{align}
Using \eqref{eqn:alphaBound}, it follows that
\begin{align}
Q^2 &\leq \rho^2 - (P_1 + P_2)^2 \leq \rho^2 - \alpha. \label{eqn:upperQQ}
\end{align}
Next, notice that the following inequality holds
\begin{align}
Q^2 &= (Q_1 + Q_2)^2 = Q_1^2 + Q_2^2 + 2Q_1Q_2 \nonumber\\
&\leq r_1^2 - P_1^2 + r_2^2 - P_2^2 + 2 \sqrt{r_1^2 - P_1^2 }\sqrt{r_2^2 - P_2^2}, \label{eqn:upperQ}
\end{align}
and, thus, there exist $Q_1$ and $Q_2$ such that $Q_1^2 \leq r_1^2 - P_1^2$ and $Q_2^2 \leq r_2^2 - P_2^2$ so that the  inequality \eqref{eqn:upperQ} holds; consequently, $\bx_i \in \cX(\underline{p}_i, \overline{p}_i, r_i)$, $i = 1, 2$.. Combining \eqref{eqn:upperQQ} and \eqref{eqn:upperQ}, we require that
$\rho^2 - \alpha \leq r_1^2 - P_1^2 + r_2^2 - P_2^2 + 2 \sqrt{r_1^2 - P_1^2 }\sqrt{r_2^2 - P_2^2},$
which, based on~\eqref{eqn:betaBound}, can be satisfied by requiring~\eqref{eq:rho} 
Thus, if $\rho$ satisfies \eqref{eq:rho}, there exists $\bx_i = [P_i, Q_i]^\sfT \in \cX(\underline{p}_i, \overline{p}_i, r_i)$, $i = 1, 2$, such that $\bar{\bx} = \bx_1 + \bx_2$, which completes the proof.

\vspace{-0.3cm}
\rev{
\subsection{Illustration of Assumption \ref{ass:dual_variable}}
\label{sec:ass2_illustration}
We next provide a simple two-dimensional example in which Assumption \ref{ass:dual_variable} is satisfied; for more elaborate examples and further conditions under which this assumption is satisfied, see \cite{Yu18}. 

Consider an example of problem \eqref{eq:aggregation} given by:
\begin{subequations} 
\label{eqn:agg_example}
\begin{align} 
\bar{f}(\bar{x})  := \quad  &\min_{x_1, x_2}  \hspace{.2cm} x_1^2  + x_2^2 \\
&\hspace{-.3cm} \textrm{subject to:} ~  x_1, x_2 \in [0, 1]\\
&\hspace{1.3cm}  x_1 + x_2 = \bar{x}
 \end{align}
\end{subequations}
for any $\bar{x} \in [0, 2]$.
We next show that the dual function associated with \eqref{eqn:agg_example} satisfies Assumption \ref{ass:dual_variable}. The dual function is given by:
\[
d(\xi) = \min_{x_1, x_2 \in [0, 1]} \left( x_1^2 + x_2^2 + \xi(x_1 + x_2 -  \bar{x}) \right) = -\xi \bar{x} +  g(\xi)
\]
where
\[
g(\xi) := \min_{x_1, x_2 \in [0, 1]} \left( x_1^2 + x_2^2 + \xi(x_1 + x_2) \right).
\]
It is easy to verify that the explicit solution of this optimization problem is given by:
\[
g(\xi) = 
\begin{cases}
0, & \textrm{ if } \xi > 0, \\
-\frac{\xi^2}{2}, & \textrm{ if } \xi \in [-2, 0] \\
2 \xi + 2, & \textrm{ if } \xi < -2,
\end{cases}
\]
and the dual function is thus given by:
\[
d(\xi) = 
\begin{cases}
-\xi \bar{x}, & \textrm{ if } \xi > 0, \\
-\frac{\xi^2}{2} -\xi \bar{x}, & \textrm{ if } \xi \in [-2, 0] \\
(2 - \bar{x}) \xi + 2, & \textrm{ if } \xi < -2.
\end{cases}
\]
Note that, for any $\bar{x} \in (0, 2)$ (i.e., in the interior of the Minkowski sum $[0, 1] + [0, 1] = [0, 2]$), the unique optimal dual variable is given by $\xi = -\bar{x}$ and it lies in a locally strongly concave region of both $d(\xi)$ and $g(\xi)$. Also, note that: (i) the dual function is not strongly concave globally, hence the standard results from, e.g.~\cite[Proposition 12.60]{rockafellar2009variational}, cannot be applied directly; and (ii) if $\bar{x} = 2$ or $\bar{x} = 0$ (i.e., a point on the boundary of the Minkowski sum), in fact there is infinite number of optimal solutions to $\max_{\xi} d(\xi)$.

}
\vspace{-0.3cm}

\subsection{Proof of Lemma \ref{lemma:strong_dual}}
\label{sec:proofthm_strong_dual}

Notice first that constraint~\eqref{eqn:constr_sum} can be rewritten in a compact form as $\bH_j \bx - \bar{\bx}_j = \mathbf{0}_2$, where $\bx \in \mathbb{R}^{2|\bar{\cD}_j|}$ stacks the vectors $\{\bx_j, j \in \bar{\cD}_j\}$ and the $2 \times 2|\bar{\cD}_j|$ matrix $\bH_j := [\bI_2, \ldots,  \bI_2]$ is full row rank.

Let $F^{(k)}(\bx) := \sum_{i \in \bar{\cD}_j} f_i^{(k)}(\bx_i)$ for brevity. From the first-order optimality conditions it follows that 
\begin{align} 
\nabla_{\bx} F|_{\bx^{\textrm{opt}}} + \bH_j^\sfT \bxi^{\textrm{opt}} = \mathbf{0} \label{eq:opt_sol}
\end{align}
where $\{\bx_i^{\textrm{opt}} \in \cX_{i}^{(k)}\}_{i \in \bar{\cD}_j}$ and $\xi^{\textrm{opt}}$ are the optimal primal and dual variables, respectively. Notice that $\bH_j^\sfT$ is a tall matrix with full column rank; therefore, its left Moore-Penrose pseudoinverse $(\bH_j^\sfT)^+$ exists. Condition~\eqref{eq:opt_sol} can thus be rewritten as $\bxi^{\textrm{opt}} = -(\bH_j^\sfT)^+ \nabla_{\bx} F|_{\bx^{\textrm{opt}}}$. Taking the norm on both sides and using the Cauchy-Schwarz inequality, one has that 
\begin{align} 
\|\bxi^{\textrm{opt}}\|_2 \leq \|(\bH_j^\sfT)^+\|_2 \|\nabla_{\bx} F|_{\bx^{\textrm{opt}}}\|_2.
\end{align}
Clearly, $\|(\bH_j^\sfT)^+\|_2 < \infty$ by construction. Also, note that: (i) $F$ is continuously differentiable by Assumption~\ref{ass:agg_cost}, and (ii) the gradient map $\nabla_{\bx} F$ is defined over a compact set. Therefore, $\nabla_{\bx} F$ is a continuous function defined over a compact set, and hence $\|\nabla_{\bx} F|_{\bx^{\textrm{opt}}}\|_2$ is bounded. This implies that  $\|\bxi^{\textrm{opt}}\|_2^2 < \infty$ as required.  Uniqueness of the optimal dual variable is implied by Assumption~\ref{ass:dual_variable}.

\vspace{-0.3cm}

\subsection{Proof of Theorem \ref{prop:derivative}}
\label{sec:proofthm_derivative}

Convexity of the optimal value function $\bar{f}_j^{(k)}(\bar{\bx})$ follows~\cite[Lemma 4.24]{RuszczynskiBook}. On the other hand,~\cite[Theorem 4.26]{RuszczynskiBook} shows that the sub-differential of $\bar{f}_j^{(k)}$ coincides with the set of optimal dual variables  associated with constraint~\eqref{eqn:constr_sum}; given~\cite[Theorem 4.26]{RuszczynskiBook} and Assumption~\ref{ass:dual_variable},~\eqref{eq:derivative_agg} then follows. 

To show that $\bar{f}_j^{(k)}(\bar{\bx})$ is Lipschitz continuous, notice first that from the convexity of the optimal value function one has that 
$\bar{f}_j^{(k)}(\bar{\bx}) - \bar{f}_j^{(k)}(\bar{\bx}^\prime) \leq (\nabla \bar{f}_j^{(k)}(\bar{\bx}))^\sfT (\bar{\bx} - \bar{\bx}^\prime)$ for any $\bar{\bx}, \bar{\bx}^\prime \in \oplus_{i \in \bar{\cD}_j} \cX_j$. 
It then follows that
\begin{subequations}
\begin{align} 
\left|\bar{f}_j^{(k)}(\bar{\bx}) - \bar{f}_j^{(k)}(\bar{\bx}^\prime)\right| & \leq \left\|\nabla \bar{f}_j^{(k)}(\bar{\bx})\right\|_2 \left\|\bar{\bx}- \bar{\bx}^\prime \right\|_2 \\
& = \| \bxi_j^{(k)} \|_2 \left\|\bar{\bx}- \bar{\bx}^\prime \right\|_2 \, .
\end{align}
\end{subequations}
Since $\bxi_j^{(k)} < \infty$ from Lemma~\ref{lemma:strong_dual}, the result readily follows.

\subsection{Proof of Theorem \ref{prop:derivative_lip}}
\label{sec:proofthm_derivative_lip}

Consider writing the dual function as
\begin{align} 
d_j(\bxi) = - \bxi^\sfT \bar{\bx}_j + g(\bxi) \label{eq:dual_rewritten}
\end{align}
where the function $g(\bxi)$ is defined as
\begin{align} 
g(\bxi) := \inf_{\{\bx_i \in \cX_i\}_{i \in \bar{\cD}_j}} F^{(k)}(\bx) + \bxi^\sfT \bH_j \bx
\end{align}
and the superscript $^{(k)}$ is dropped for brevity. From Assumption~\ref{ass:dual_variable} (see also~\cite{Yu18}), it follows that $g(\bxi)$ is locally strongly concave and differentiable; denote as $\beta > 0$ the (local) strong concavity coefficient.  

For any feasible $\bar{\bx}_j$ and $\bar{\bx}_j^\prime$ in the interior of the Minkowski sum of $\cX_{i}, \,  i \in \bar{\cD}_j$, let $\bxi^{\star}$ and $\bxi^{\prime \star}$ denote the corresponding optimal dual variables. From the optimality of $\bxi^{\star}$ and $\bxi^{\prime \star}$, we have that
\begin{align} 
(\nabla_\bxi g(\bxi^{\star}) - \bar{\bx}_j )^\sfT (\bxi - \bxi^{\star}) \leq 0 , \hspace{.2cm} \forall \, \bxi \, . \label{eq:opt_xi} \\
(\nabla_\bxi g(\bxi^{\prime \star}) - \bar{\bx}_j^\prime )^\sfT (\bxi - \bxi^{\prime \star}) \leq 0 , \hspace{.2cm} \forall \, \bxi \, . \label{eq:opt_xi_prime}
\end{align}
By using $\bxi = \bxi^{\prime \star}$ in \eqref{eq:opt_xi} and $\bxi = \bxi^{\star}$ in  \eqref{eq:opt_xi_prime}, and summing up these two inequalities, we obtain
\begin{subequations}
\begin{align} 
(\bar{\bx}_j - \bar{\bx}_j^\prime)^\sfT (\bxi^{\star} - \bxi^{\prime \star}) & \leq \left(\nabla_\bxi g(\bxi^{\star}) - \nabla_\bxi g(\bxi^{\prime \star}) \right)^\sfT (\bxi^{\star} - \bxi^{\prime \star}) 
\label{eq:proof_4_1} \\
& \leq - \beta \| \bxi^{\star} - \bxi^{\prime \star} \|_2^2 \label{eq:proof_4_2},
\end{align}
\end{subequations}
where the last inequality follows by the local  strong concavity of $g(\bxi)$ around the optimal dual variables. This implies
\begin{subequations}
\begin{align} 
 \| \bxi^{\star} - \bxi^{\prime \star} \|_2^2 & \leq 
 \frac{1}{\beta}| (\bar{\bx}_j - \bar{\bx}_j^\prime)^\sfT (\bxi^{\star} - \bxi^{\prime \star}) | \\
 & \leq \frac{1}{\beta}\|\bar{\bx}_j - \bar{\bx}_j^\prime \|_2 \| \bxi^{\star} - \bxi^{\prime \star}\|_2
\end{align}
\end{subequations}
where the second inequality follows from the Cauchy-Schwarz inequality. Therefore, whenever $\bxi^{\star} \neq \bxi^{\prime \star}$, we have that
\begin{align} 
\| \bxi^{\star} - \bxi^{\prime \star} \|_2 \leq \frac{1}{\beta} \|\bar{\bx}_j - \bar{\bx}_j^\prime \|_2
\end{align}
which proves the theorem.

\subsection{Proof of Theorem \ref{theorem.inexact}}
\label{sec:proofthm_inexact}

The proof of Theorem \ref{theorem.inexact} follows steps that are similar to the ones outlined in~\cite{opfPursuit}. Consider the map $\bPhi^{(k)}$ defined as:
\begin{equation*}
\label{eq:phi_mapping}
 \bPhi^{(k)}: \{\bz^{t_k}\} \mapsto 
 \left[\begin{array}{c}
 \nabla_{[\bx, \bar{\bx}]}  L_r^{(k)}(\bx, \bar{\bx}, \bd) |_{\bx^{(k)}, \bar{\bx}^{(k)}, \bd^{(k)}} \\
- \nabla_{\bd}  L_r^{(k)}(\bx, \bar{\bx}, \bd) |_{\bx^{(k)}, \bar{\bx}^{(k)}, \bd^{(k)}} 
\end{array}
\right],
\end{equation*}
and notice from~\cite{opfPursuit,Koshal11} that $\bPhi^{(k)}$ is strongly monotone with constant $\min\{r_p,r_d\}$, and Lipschitz over the domain of the primal and dual variables with constant $L_\Phi = [(L+r_p+5G)^2 + 5(G+r_d^2]^{\frac{1}{2}}$. Next, let $\bPhi_e^{(k)}$ denote the counterpart when feedback is utilized in the gradient computation, and consider the following  inequality: 
\begin{align}
& \hspace{-.3cm} \|\bz^{(k)} - \bz^{(k-1,\star)} \|_2 \leq \| \hat{\bz}^{(k-1)}  - \alpha \bPhi_e^{(k-1)}(\hat{\bz}^{(k-1)}) \nonumber \\
& \hspace{2.3cm} - \bz^{(k-1,\star)} +  \alpha \bPhi^{(k-1)}(\bz^{(k-1,\star)})  \|_2 \, . \label{eq:proof_thm5_1}
\end{align}
Recognizing that $\bPhi_e^{(k)}(\hat{\bz}^{(k)}) - \bPhi^{(k)}(\bz^{(k)}) = \be^{(k)}$, and adding and subtracting $\bz^{(k-1)}$ on the right-hand-side of~\eqref{eq:proof_thm5_1}, it follows that~\eqref{eq:proof_thm5_1} can be further bounded as:
\begin{align}
& \|\bz^{(k)} - \bz^{(k-1,\star)} \|_2 \leq \|\hat{\bz}^{(k-1)} - \bz^{(k-1)}\|_2 +  \alpha \|\be^{(k-1)}\|_2 \nonumber \\
& + \|\bz^{(k-1)}  - \alpha \bPhi^{(k-1)}(\bz^{(k-1)}) - \bz^{(k-1,\star)} +  \alpha \bPhi^{(k)}(\bz^{(k-1,\star)})  \|_2 \, . \label{eq:proof_thm5_2}
\end{align}
Following~\cite{opfPursuit}, the third term on the right-hand-side of~\eqref{eq:proof_thm5_2} can be bounded with the term $c(\alpha) \|\bz^{(k-1)} - \bz^{(k-1),\star}\|_2$; hence,
\begin{align}
& \|\bz^{(k)} - \bz^{(k-1,\star)} \|_2 \leq e_x^{(k-1)} + \alpha \|\be^{(k-1)}\|_2 \nonumber \\
& \hspace{1.5cm} + c(\alpha) \|\bz^{(k-1)} - \bz^{(k-1),\star}\|_2 . \label{eq:proof_thm5_3}
\end{align}
Next, $\|\bz^{(k)} - \bz^{(k,\star)}\|_2$ can be bounded as:
\begin{align}
\|\bz^{(k)} - \bz^{(k,\star)}\|_2 & = \|\bz^{(k)} - \bz^{(k-1,\star)} + \bz^{(k-1,\star)} - \bz^{(k,\star)} \|_2 \nonumber \\
& \hspace{-2.5cm} \leq \|\bz^{(k-1,\star)} - \bz^{(k,\star)} \|_2  + \|\bz^{(k)} - \bz^{(k-1,\star)}\|_2  \label{eq:proof_thm5_4} \\
& \hspace{-2.5cm}  \leq \sigma^{(k)} + e_x^{(k-1)} + \alpha \|\be^{(k-1)}\|_2 + c(\alpha) \|\bz^{(k-1)} - \bz^{(k-1),\star}\|_2  \label{eq:proof_thm5_5} .
\end{align}
By recursively applying~\eqref{eq:proof_thm5_5}, the result of Theorem \ref{theorem.inexact} follows.

\subsection{Data for the Simulations}
\label{sec:data}

Table~\ref{tab:pv_systems} and Table~\ref{tab:battery_systems} list the locations, capacities, and connection types of PV systems and energy storage systems in the considered distribution network.

 \begin{table}[t]
\begin{center}
\caption{PV systems in the feeder}
\label{tab:pv_systems}
\begin{tabular}{|c|r|c|c||c|c|c|c|}
\hline 
\textbf{Node} & \textbf{kW} & \textbf{Con.} & \textbf{Num.} & \textbf{Node} & \textbf{kW} & \textbf{Con.} & \textbf{Num.}  \\
\hline \hline 
8 & 61.57 & 3$\Delta$ & 5 & 82 & 154.00 & 3$\Delta$ & 2 \\
9 & 71.95 & 3$\Delta$ & 4 & 85 & 121.59 & 3$\Delta$ & 1 \\
10 & 600.00 & 3$\Delta$ & 2 & 86 & 188.13 & 3$\Delta$ & 2 \\
16 & 11.53 & 3$\Delta$ & 1 & 88 & 19.15 & 3$\Delta$ & 2 \\
19 & 66.05 & 3$\Delta$ & 1 & 100 & 36.94 & 1$\Delta$ & 3 \\
21 & 18.27 & 3$\Delta$ & 2 & 93 & 12.95 & 3$\Delta$ & 1 \\
23 & 13.53 & 3$\Delta$ & 1 & 99 & 116.27 & 3$\Delta$ & 2 \\
32 & 100.00 & 3$\Delta$ & 1 & 100 & 13.06 & 3$\Delta$ & 1 \\
40 & 6.35 & 3$\Delta$ & 1 & 105 & 700.00 & 3$\Delta$ & 1 \\
44 & 106.14 & 3$\Delta$ & 1 & 110 & 31.17 & 3$\Delta$ & 1 \\
48 & 293.54 & 3$\Delta$ & 1 & 112 & 22.63 & 3$\Delta$ & 3 \\
58 & 2.20 & 1Y & 1 &  115 & 691.07 & 3$\Delta$ & 3 \\
63 & 2.20 & 1$\Delta$ & 1 & 117 & 9.76 & 3$\Delta$ & 2 \\
65 & 9.97 & 3$\Delta$ & 1 & 119 & 7.66 & 3$\Delta$ & 1 \\
67 & 5.7 & 1$\Delta$ & 1 & 121 & 100.00 & 3$\Delta$ & 1 \\
68 & 4.6 & 1$\Delta$ & 1 & 123 & 19.35 & 3$\Delta$ & 1 \\
& & & & 125 & 50.23 & 3$\Delta$ & 3 \\  
\hline 
\end{tabular}
\end{center}
\end{table}
 \begin{table}[t]
\begin{center}
\caption{Batteries in the feeder}
\label{tab:battery_systems}
\begin{tabular}{|c|r|r|c|c|}
\hline 
\textbf{Node} & \textbf{kW} & \textbf{kWh} & \textbf{Con.} & \textbf{Num.} \\
\hline \hline 
8 & 22.98 & 35.05 & 3$\Delta$ & 5 \\
9 & 13.73 & 7.17 & 3$\Delta$ & 4 \\
10 & 46.85 & 39.53 & 3$\Delta$ & 1 \\
16 & 3.40 & 3.07 & 3$\Delta$ & 1 \\
19 & 26.11 & 74.68 & 3$\Delta$ & 1 \\
21 & 13.79 & 21.46 & 3$\Delta$ & 2 \\
23 & 4.60 & 19.68 & 3$\Delta$ & 1 \\
32 & 60.61 & 60.23 & 3$\Delta$ & 1 \\
40 & 3.04 & 0.82 & 3$\Delta$ & 1 \\
44 & 4.80 & 3.01 & 3$\Delta$ & 1 \\
48 & 49.77 & 17.75 & 3$\Delta$ & 1 \\
65 & 11.89 & 3.20 & 3$\Delta$ & 1 \\
82 & 52.74 & 73.81 & 3$\Delta$ & 1 \\
85 & 5.44 & 6.37 & 3$\Delta$ & 1 \\
86 & 7.41 & 5.19 & 3$\Delta$ & 2 \\
88 & 7.19 & 5.19 & 3$\Delta$ & 2 \\
90 & 10.39 & 13.85 & 3$\Delta$ & 1 \\
93 & 3.23 & 2.40 & 3$\Delta$ & 1 \\
99 & 13.84 & 12.23 & 3$\Delta$ & 2 \\
100 & 6.40 & 5.58 & 3$\Delta$ & 1 \\
110 & 3.02 & 0.99 & 3$\Delta$ & 1 \\
112 & 15.13 & 7.82 & 3$\Delta$ & 3 \\
115 & 31.17 & 89.18 & 3$\Delta$ & 3 \\
117 & 6.63 & 5.32 & 3$\Delta$ & 2 \\
119 & 5.75 & 1.54 & 3$\Delta$ & 1 \\
123 & 5.09 & 1.92 & 3$\Delta$ & 1 \\
125 & 13.83 & 14.51 & 3$\Delta$ & 3 \\
\hline 
\end{tabular}
\end{center}
\end{table}

\bibliographystyle{IEEEtran}
\bibliography{biblio.bib}

\end{document}